\documentclass{article}
\usepackage{latexsym,amssymb}

\newcommand{\ta}{{\tilde a}}
\newcommand{\Nabla}{{\bigtriangledown}}

\newcommand{\kerr}{{\rm ker}}
\newcommand{\im}{{\rm im}}

\newcommand{\tg}{{\Tilde{g}}}

\newcommand{\hg}{{\widehat{g}}}

\newcommand{\ev}{{{\rm ev}}}
\newcommand{\SFr}{{\rm SFr}}
\newcommand{\SO}{{\rm SO}}
\newcommand{\R}{{\mathbb R}}
\newcommand{\RP}{{\mathbb RP}}
\newcommand{\T}{{\mathbb T}}

\newcommand{\Q}{{\mathbb Q}}

\newcommand{\Z}{{\mathbb Z}}
\newcommand{\C}{{\mathbb C}}
\newcommand{\CP}{{\mathbb CP}}
\newcommand{\U}{{\mathbb U}}

\newcommand{\PU}{{\mathbb PU}}

\newcommand{\Mm}{{\mathcal M}}
\newcommand{\Hh}{{\mathcal H}}

\newcommand{\Ss}{{\mathcal S}}
\newcommand{\Ll}{{\mathcal L}}
\newcommand{\Uu}{{\mathcal U}}
\newcommand{\PPp}{{\mathcal P}}

\newcommand{\oMm}{{\overline{\mathcal M}}}
\newcommand{\oMmn}{{\overline{\mathcal M}\,\!^\nu}}
\newcommand{\pJ}{{\overline{\partial}_J}}
\newcommand{\ov}{\overline}

\newcommand{\bcp}{{\ov{\CP}\,\!^2}}

\newcommand{\Gg}{{\mathcal G}}

\newcommand{\Tilde}{\widetilde}
\newcommand{{\TS}}{{\Tilde {\rm Symp_0}}}
\newcommand{{\TIH}}{{\Tilde {\rm Ham}}}
\newcommand{\p}{{\partial}}
\newcommand{\al}{{\alpha}}

\newcommand{\be}{{\beta}}

\newcommand{\Om}{{\Omega}}
\newcommand{\om}{{\omega}}
\newcommand{\eps}{{\varepsilon}}
\newcommand{\de}{{\delta}}

\newcommand{\ga}{{\gamma}}
\newcommand{\Ga}{{\Gamma}}
\newcommand{\io}{{\iota}}
\newcommand{\ka}{{\kappa}}
\newcommand{\la}{{\lambda}}

\newcommand{\si}{{\sigma}}

\newcommand{\Symp}{{\rm Symp}}
\newcommand{\Diff}{{\rm Diff}}
\newcommand{\Ham}{{\rm Ham}}

\newcommand{\LHam}{{\rm LHam}}
\newcommand{\Hom}{{\rm Hom}}
\newcommand{\Flux}{{\rm Flux}}
\newcommand{\Hor}{{\it Hor}}
\newcommand{\Aut}{{\rm Aut}}
\newcommand{\tr}{{\rm tr}}

\newcommand{\Si}{{\Sigma}}

\newcommand{\MS}{{\medskip}}
\newcommand{\SmS}{{\smallskip}}

\newcommand{\NI}{{\noindent}}
\newcommand{\proof}[1]{\noindent{\bf Proof#1:\  }}
\newcommand{\QED}{\hfill$\Box$\medskip}
\newtheorem{theorem}{Theorem}[section]
\newtheorem{thm}[theorem]{Theorem}
\newtheorem{prop}[theorem]{Proposition}
\newtheorem{cor}[theorem]{Corollary}
\newtheorem{lemma}[theorem]{Lemma}

\newtheorem{defn}[theorem]{Definition}

\newtheorem{rmk}[theorem]{Remark}

\newtheorem{ex}[theorem]{Example}

\begin{document}

\title{Symplectic structures on fiber bundles}
\author{Fran\c{c}ois Lalonde\thanks{Partially supported by NSERC grant 
OGP 0092913
and FCAR grant ER-1199.} \\ Universit\'e du Qu\'ebec \`a Montr\'eal
\\ (flalonde@math.uqam.ca) \and Dusa McDuff\thanks{Partially
supported by NSF grants DMS 9704825 and DMS 0072512.} \\ State University of New York
at Stony Brook \\ (dusa@math.sunysb.edu)
}

\date{Aug 27, 2000, revised March 9, 2001, with erratum Dec 2004}
\maketitle

\MS
\MS
\MS

\NI
{\bf Abstract}

 Let $\pi: P\to B$ be a locally trivial fiber bundle 
over a connected CW complex $B$
with fiber equal to the closed symplectic manifold 
$(M,\om)$. Then $\pi$  is said to be a
{\it symplectic  fiber bundle} if its structural group  is
the group of  symplectomorphisms $\Symp(M,\om)$,  
and is called {\it Hamiltonian} if this group may be reduced to
the group $\Ham(M,\om)$ of Hamiltonian symplectomorphisms.
In this paper, building on prior work by Seidel and Lalonde, McDuff and
Polterovich, we show that these bundles have interesting cohomological 
properties. In particular,  for many bases $B$ (for example when $B$ is a
sphere, a coadjoint orbit or a product of complex projective spaces) 
the rational cohomology
of $P$ is the tensor product of the cohomology of $B$ with that of $M$.  
As a consequence the natural action of the rational  homology 
$H_k(\Ham(M))$ on $H_*(M)$ is trivial for all $M$ and all $k > 0$.

Added:  The erratum makes a small change to Theorem 1.1 that characterizes Hamiltonian bundles.\MS\MS
\MS
\MS
\MS

\NI
keywords: symplectic fiber bundle, Hamiltonian fiber bundle, symplectomorphism
group, group of Hamiltonian symplectomorphisms, rational cohomology of fiber bundles

\MS

\NI

\newpage

\tableofcontents

\section{Introduction and main results}\label{se:intro}

In this section we first discuss how to characterize
 Hamiltonian bundles and their automorphisms, and then
describe their main properties, in particular deriving
conditions under which 
the cohomology of the total space splits as a
product.   Finally we state some applications to
the action of $\Ham(M)$ on $M$ and to nonHamiltonian 
symplectic bundles.
This paper should be considered as a sequel to
Lalonde--McDuff--Polterovich~\cite{LMP} and McDuff~\cite{Mc} 
which establish analogous results for 
Hamiltonian bundles over $S^2$.
Several of our results are well known for Hamiltonian bundles
whose structural group is a compact Lie group. They therefore fit in with the
idea mentioned by Reznikov~\cite{Rez}
 that the group of symplectomorphisms behaves cohomologically much like a
Lie group.  \MS

\subsection{Characterizing Hamiltonian bundles}\label{ss:char}

A fiber bundle $M\to P\to B$ is said to be symplectic if its structural group
reduces to the group of symplectomorphisms $\Symp(M, \om)$ of 
the closed symplectic manifold $(M,\om)$.
In this case, each fiber $M_b = \pi^{-1}(b)$ is equipped with a well defined
symplectic form $\om_b$ such that $(M_b,\om_b)$ is symplectomorphic to 
$(M, \om)$.  Our
first group of results establish geometric criteria for a symplectic bundle
to be Hamiltonian, i.e. for the structural group to reduce to $\Ham(M, \om)$. 
Quite often we simplify the notation by writing $\Ham(M)$ and $\Symp_0(M)$ 
(or even $\Ham$ and $\Symp_0$)
instead of $\Ham(M, \om)$ and $\Symp_0(M, \om)$.

Recall that the group  $\Ham(M, \om)$ is a
connected normal subgroup of the identity
component $\Symp_0(M,\om)$ of the group of 
symplectomorphisms, and fits into the exact sequence
$$
\{id\} \to \Ham(M, \om) \to \Symp_0(M, \om)\stackrel{\Flux}\to 
H^1(M,\R)/\Ga_{\om}  \to
\{0\},
$$
where $\Ga_{\om}$ is the flux group.\footnote
{
It is not known whether this group is discrete in all cases,
although its rank is always less
than or equal to the first  Betti number of $M$ 
by the results of~\cite{LMP}.
It is  discrete if $[\om]$ is a 
rational (or integral) class and in various other cases discussed 
in~\cite{LMP2} and Kedra~\cite{Ke}.} 
Because $\Ham(M)$ is connected,
every Hamiltonian bundle is symplectically 
trivial over the $1$-skeleton of the base.    The
following proposition was proved in~\cite{MS}~Thm.~6.36 by a somewhat
analytic argument. We give a more topological proof in~\S\ref{ss:hamchar}
below.

\begin{thm}\label{thm:hamchar}  A symplectic bundle
  $\pi: P\to B$ is 
Hamiltonian
 if and only if the
following conditions hold:\begin{itemize}
 \item[(i)]  the restriction of $\pi$ to the $1$-skeleton $B_1$ of $B$
is  symplectically trivial, and
\item[(ii)] there is
a cohomology class $a\in H^2(P,\R)$ that restricts to $[\om_b]$ on $M_b$.
\end{itemize}
\end{thm}

There is no loss of generality in assuming that the bundle
$\pi: P\to B$ is smooth. Then
recall from 
Guillemin--Lerman--Sternberg~\cite{GLS} (or~\cite{MS} Chapter 6) that any 
$2$-form $\tau$ on $P$  that restricts to $\om_b$ on each fiber $M_b$ defines 
 a connection $\Nabla_\tau$  on $P$ whose horizontal distribution $\Hor_\tau$
is just the $\tau$-orthogonal complement of the tangent spaces
to the fibers:
$$
\Hor_\tau(x) = \{v\in T_xP: \tau(v,w) = 0\;\mbox{for all } w\in
T_x(M_{\pi(x)})\}. $$
Such forms $\tau$ are called {\it connection forms}.
The closedness of $\tau$ is a sufficient (but not
necessary) condition for the holonomy  of $\Nabla_\tau$ to be
symplectic, see Lemma~\ref{le:conn}.  A simple argument due to Thurston
(\cite{MS}~Thm. 6.3 for instance)  shows that the 
cohomological condition $(ii)$ above is equivalent to the existence of a 
closed extension $\tau$ of the forms $\om_b$.  Condition $(i)$ is then
equivalent to requiring that the holonomy 
of $\Nabla_\tau$
around any loop in $B$ 
  belongs to the identity component $\Symp_0(M)$ 
of $\Symp(M)$. 
 Hence  the above result can be  rephrased in terms of such closed extensions
$\tau$ as follows.

\begin{prop}\label{prop:hamchar2}
A  symplectic
bundle  $\pi: P\to B$ is Hamiltonian if and only if the forms $\om_b$ 
on the fibers have  a
closed extension $\tau$ such that the  holonomy of $\Nabla_\tau$
around any loop in $B$ 
lies in the identity component $\Symp_0(M)$ of $\Symp(M)$. 
\end{prop}

This is a slight extension of a result of 
Guillemin--Lerman--Sternberg, who called a specific choice of
$\tau$ the {\it coupling form}: see also~\cite{MS}~Thm.~6.21.
As we show in \S\ref{ss:funct} the existence of $\tau$ is the key to the good
behavior of Hamiltonian bundles under 
composition.

\begin{rmk}\rm\label{rmk:sconn}
When $M$ is simply connected, $\Ham(M)$ is the
identity component $\Symp_0(M)$ of $\Symp(M)$, and so a symplectic bundle is
Hamiltonian if and only if condition $(i)$ above is satisfied, i.e. if and only
if it is trivial over the $1$-skeleton $B_1$. In this case,
as observed by Gotay {\it et al.} in~\cite{GLSW}  Theorem~2,
it is known that $(i)$ implies $(ii)$ 
for general topological reasons to do with
the behavior of evaluation maps.  (One can reconstruct their arguments
from our  Lemmas~\ref{le:shamchar} and~\ref{le:ev2}.)
More generally, $(i)$ implies $(ii)$  for all symplectic bundles 
 with fiber $(M,\om)$ if
and only if the flux group  $\Ga_{\om}$ vanishes. \end{rmk}

\subsection{Hamiltonian structures and their
automorphisms}\label{ss:strchar}

The question then arises as to what a Hamiltonian structure on a fiber bundle
actually is. How many Hamiltonian structures can one put on a given 
symplectic bundle $\pi:P\to B$?
What does one mean by an automorphism of such a structure? 
These questions are discussed in detail  in~\S\ref{ss:hamstr}.
We now  summarize the results of that discussion.

In homotopy theoretic terms, a Hamiltonian
structure on a symplectic bundle  $\pi:P\to B$
 is simply a lift $\tg$ to $B\,\Ham(M)$ of
the classifying map $g: B \to B\,\Symp(M, \om)$ of the underlying symplectic
bundle, i.e.  it is a
homotopy commutative diagram
$$
\begin{array}{ccc}
& & B\,\Ham(M)\\
& { \tg} \nearrow& \downarrow\\
B &\stackrel{g}\longrightarrow & B\,\Symp(M).
\end{array}$$
Hamiltonian structures are in bijective correspondence 
with homotopy classes of such lifts.  
There are two stages to choosing the lift $\tg$:  one first  lifts $g$ to
a map   $\hg$ into $B\,\Symp_0(M,\om)$, where $\Symp_0$ is the identity
component of $\Symp$, and then to a map $\tg$ into $B\,\Ham(M,\om)$.  
As we show in~\S\ref{ss:hamstr}, choosing  $\hg$ is equivalent to fixing 
the isotopy class of an
identification of $(M, \om )$ with the fiber $(M_{b_0}, \om_{b_0})$ over the
base point $b_0$. If $B$ is simply connected, in particular if $B$ is a single
point,  there is then a unique Hamiltonian
structure on $P$, i.e. a unique choice of lift $\tg$.   Before describing what
happens in the general case, we discuss properties of the
extensions $\tau$.
 
 Let $\tau \in \Om^2(P)$ be
a closed extension of the symplectic forms on the fibers. 
Given a loop $\ga: S^1\to  B$ based at 
$b_0$, and  a symplectic trivialisation $T_\ga: \ga^*(P) \to S^1 \times (M,
\om)$ that extends the given identification of $M_{b_0}$ with $M$,
push forward $\tau$ to a form $(T_\ga)_*\tau$ 
on $S^1 \times (M, \om)$. Its
characteristic flow round $S^1$ is transverse to the fibers and defines
a symplectic isotopy $\phi_t$ of $(M, \om ) = (M_{b_0}, \om_{b_0})$
whose flux, as map from $H_1(M) \to \R$, is equal to 
$(T_\ga)_*[\tau] ([S^1] \otimes \cdot)$: see Lemma~\ref{le:conn2}. 
This flux depends only on the cohomology class $a$ of $\tau$.
Moreover, as we mentioned above, any extension $a$ of the fiber class
$[\om]$ can be represented by a form $\tau$ that extends the $\om_b$.
 Thus, given $T_\ga$ and an extension $a = [\tau] \in H^2(P)$ of the
fiber symplectic class $[\om]$, it makes sense to define
the {\it flux class}  $f(T_\ga, a) \in H^1(M,\R)$ by 
$$
f(T_\ga, a)(\de) = (T_\ga)_*(a)(\ga \otimes \de) \; \; \; {\rm for} \; {\rm all} \; 
\de \in H_1(M).
$$
The equivalence class $[f(T_\ga,a)] \in H^1(M,\R)/\Ga_{\om}$
does not depend on the choice of $T_\ga$: indeed two such choices differ
by a loop $\phi$ in $\Symp_0(M,\om)$ and so the difference
$$
f(T_\ga, a) - f(T_\ga', a) = f(T_\ga,a) \circ \tr_{\phi} = \om \circ \tr_{\phi}
$$
belongs to $\Ga_{\om}$. 
 The following lemma is elementary:
see~\S\ref{ss:hamstr} : 

\begin{lemma} \label{le:normalized} If $\pi:P\to B$ is a symplectic 
bundle satisfying the
conditions of Theorem~\ref{thm:hamchar}, there is an extension $a$ of the
symplectic fiber class that has
trivial flux 
$$
[f(T_\ga,a)] = 0\;\; \in\;\; H^1(M,\R)/\Ga_{\om}
$$
round each loop $\ga$ in $B$.
\end{lemma}

\begin{defn}\label{def:norm}\rm  An extension $a$ of the
symplectic fiber class  $[\om_{b_0}]$  is  {\it normalized} if it satisfies the
conclusions of the above lemma.  Two such extensions $a$ and $a'$ are
equivalent  (in symbols, $a\sim a'$)   if and only if they
have equal restrictions to $\pi^{-1}(B_1)$, or, equivalently, if and only if
$a - a' \in \pi^*(H^2(B))$.
\end{defn}

We show in~\S\ref{ss:hamstr} that
Hamiltonian structures are
in one-to-one correspondence with symplectic trivializations of the
$1$-skeleton $B_1$ of $B$, with two such trivializations being equivalent if
and only if they differ by Hamiltonian loops.   
If two trivializations
$T_\ga, T'_\ga$ differ by a Hamiltonian loop $\phi$ then $f(T_\ga, a) -
f(T_\ga', a) = 0$.   In terms of fluxes of closed extensions, we therefore get:

\begin{thm}\label{thm:hamchar3}
Assume that  a symplectic bundle $\pi:P\to B$ can be 
symplectically trivialized over $B_1$. Then a Hamiltonian structure exists
on $P$ if and only if there is a normalized
extension $a$ of $\om$. Such a structure
consists of an isotopy class of  symplectomorphisms $(M, \om) \to (M_{b_0},
\om_{b_0})$ together with an equivalence class 
$\{a\}$ of normalized extensions of
the fiber symplectic class.  \end{thm}

In other words, with respect to a fixed trivialization over $B_1$, 
Hamiltonian structures are in one-to-one correspondence with
homomorphisms  $\pi_1(B) \to \Ga_{\om}$ , given by the fluxes
$f_\ga(T, a)$ of monodromies round the loops of the base. We will call
$\{a\}$ the {\it Hamiltonian extension class}, and will denote the
Hamiltonian structure on $P$ by the triple $(P, \pi, \{a\})$.  
A different description of a
Hamiltonian structure is sketched in Appendix A.
\MS

We now turn to the question of describing automorphisms of Hamiltonian
structures. It is convenient to distinguish between symplectic and Hamiltonian
automorphisms, just as we distinguish between 
$\Symp(M,\om)$ and $\Ham(M,\om)$ in the case when $B = pt.$
Notice that if $P\to B$ is a symplectic bundle, there is a natural notion of 
symplectic automorphism.  This is a
fiberwise diffeomorphism $\Phi:P\to P$ that covers the identity map on the
base and restricts on each fiber to an element $\Phi_b$ of the group
$\Symp(M_b, \om_b)$.\footnote
{
One could allow more general automorphisms of the base, but we will restrict
to this simple case here.
}
Because $\Ham(M, \om)$ is a normal subgroup of 
$\Symp(M,\om)$, it also makes sense to require that $\Phi_b\in \Ham(M_b,
\om_b)$ for each $b$.  Such automorphisms are called {\it Hamiltonian
automorphisms} of the symplectic bundle $P\to B$. Let us write $\Symp(P, \pi)$
and $\Ham(P,\pi)$ for the groups of such automorphisms.  Observe that
the group $\Ham(P,\pi)$ may not be connected.
Because the fibers of  Hamiltonian
 bundles are identified with $(M, \om)$ up to isotopy,  we shall also
need to consider the (not necessarily connected)  group $\Symp_0(P,\pi)$ of
symplectomorphisms of $(P,\pi)$ where $\Phi_b\in \Symp_0(M_b, \om_b)$
for one and hence all $b$. \MS

Now let us consider automorphisms of Hamiltonian bundles.  
As a guide note that in the trivial case when  $B = pt$, a Hamiltonian structure
on $P$ is an identification of $P$ with $M$ up to symplectic isotopy.  Hence the
group of automorphisms of this structure can be identified with
$\Symp_0(M,\om)$.  In general, 
 if $\{a\}$ is a Hamiltonian structure on $(P, \pi)$ and 
$\Phi\in \Symp_0(P, \pi)$ then
$\Phi^*(\{a\}) =  \{a\}$ if and only if $\Phi^*(a) =  a$ for some $a$ in the class
$\{a\}$, because 
$\Phi$ induces the identity map on the base and $a - a' \in \pi^*(H^2(B))$ 
when $a\sim a'$.  We therefore make the
following definition.

 \begin{defn}\label{def:hamaut}\rm
Let $(P, \pi, \{a\})$ be a Hamiltonian structure on the symplectic bundle $P\to
B$ and let
$\Phi\in \Symp(P, \pi)$.  Then $\Phi$ is an {\it automorphism of the
Hamiltonian structure} $(P, \pi, \{a\})$ if  $\Phi\in \Symp_0(P, \pi)$
and $\Phi^*(\{a\}) = \{a\}$. The group formed by these elements is denoted by
$\Aut(P,\pi,\{a\})$. 
\end{defn}

The following result is  not hard to prove, but is easiest to see in 
the context of a discussion of the action of $\Ham(M)$ on $H^*(M)$.
Therefore the proof is deferred to~\S\ref{ss:act}.

\begin{prop} \label{prop:hamchar4}
Let $P\to B$ be a Hamiltonian bundle and $\Phi\in \Symp_0(P, \pi)$.
Then the following statements are equivalent:
\begin{itemize}
\item[(i)] $\Phi$ is isotopic to an element of $\Ham(P, \pi)$;
\item[(ii)]
$\Phi^*(\{a\}) = \{a\}$ for some
Hamiltonian structure $\{a\}$  on $P$;  
\item[(iii)]
$\Phi^*(\{a\}) = \{a\}$  for all  
Hamiltonian structures $\{a\}$  on $P$. 
\end{itemize} 
\end{prop}

\begin{cor}  For any Hamiltonian   bundle
$P\to B$, the group $\Aut(P,\pi,\{a\})$ does not depend
on the choice of the Hamiltonian structure $\{a\}$ put on $P$. Moreover,
it contains $\Ham(P, \pi)$ and each element
of $\Aut(P,\pi,\{a\})$ is isotopic to an element in $\Ham(P,\pi)$.
\end{cor}

The following characterization is now obvious:

\begin{lemma}\label{le:auttriv}  Let $P$ be the product $B\times M$
and $\{a\}$ any Hamiltonian
structure.  Then:
\begin{itemize}
\item[(i)]
$\Ham(P, \pi)$ consists of all maps from $B$ to $
\Ham(M, \om)$. 
\item[(ii)]  $\Aut(P, \pi, \{a\})$ consists of all maps
$\Phi: B \to  \Symp_0(M, \om)$ for which the composite
$$
\pi_1(B)\stackrel{\Phi_*}\longrightarrow 
\pi_1(\Symp_0(M))\stackrel{\Flux_\om}\longrightarrow H^1(M,\R) $$
is trivial.
\end{itemize}
\end{lemma}

The basic reason why Proposition~\ref{prop:hamchar4} holds is that Hamiltonian
automorphisms of $(P, \pi)$ act trivially on the set of extensions of the fiber
symplectic class.   This need not be true for symplectic
automorphisms.  For example, if $\pi: P = S^1\times M \to S^1$ is a trivial
bundle and  $\Phi$ is given by a nonHamiltonian loop  $\phi$ in
$\Symp_0(M)$, then $\Phi$ is in $ \Symp_0(P, \pi)$ but it preserves
no Hamiltonian structure on $P$ since
$\Phi^*(a) = a + [dt]\otimes \Flux (\phi)$.

In general, if we choose a trivialization of $P$ over $B_1$, there are exact
sequences
$$
\{id\}  \to  \Aut(P,\pi, \{a\})  \to  
\Symp_0(P,\pi)   \to  
\Hom(\pi_1(B),\Ga_{\om}) \to   \{id\}, $$
$$
\{id\} \to   \Ham(P,\pi, \{a\})  \to  
\Aut(P,\pi, \{a\}) 
 \to   H^1(M,\R)/\Ga_{\om}  \to   \{0\}.
$$
In particular,  the subgroup of $\Aut(P,\pi,\{a\})$
consisting of automorphisms that belong to $\Ham(M_{b_0}, \om_{b_0})$ at the
base point $b_0$ retracts to 
$\Ham(P,\pi,\{a\})$.
\MS

\subsection{Stability}\label{ss:stabil}

Another important property of Hamiltonian bundles is  stability.

\begin{defn}\rm 
A symplectic (resp. Hamiltonian) bundle $\pi: P\to B$ with fiber
$(M,\om)$ is said to be {\it stable} if $\pi$ may be given a 
symplectic (resp. Hamiltonian) structure with respect to any symplectic
form $\om'$ on $M$ that is sufficiently close to (but not necessarily
cohomologous to) $\om$, in such a way that the structure depends 
continuously on $\om'$.
\end{defn}

Using Moser's homotopy argument, it is easy to prove that
any symplectic bundle is stable (see Corollary~\ref{cor:stab0}). 
The following characterization of Hamiltonian stability is an almost 
immediate consequence of
Theorem~\ref{thm:hamchar}.  It is proved in \S\ref{ss:stab} below.

\begin{lemma}\label{le:stab}  A Hamiltonian bundle $\pi: P\to
B$ is stable if and only if  the restriction map $H^2(P, \R)\to H^2(M,\R)$ is
surjective. \end{lemma}

The following result is less immediate.

\begin{thm}\label{thm:stable} Every Hamiltonian bundle is stable.
\end{thm}

The proof uses the (difficult) stability property for Hamiltonian
bundles over $S^2$ that was established in~\cite{LMP,Mc} 
as well as the (easy)
fact that the image of the evaluation map 
$\pi_2(\Ham(M)) \to \pi_2(M)$ lies in
the kernel of $[\om]$: see Lemma~\ref{le:ev2}.

\subsection{Cohomological splitting}\label{ss:split}

We next  extend the splitting results of
Lalonde--McDuff--Polterovich~\cite{LMP} and McDuff~\cite{Mc}.  
These papers prove
that the rational cohomology of every Hamiltonian bundle $\pi:P\to S^2$
 splits additively, i.e. there is an additive isomorphism
$$
H^*(P)\cong H^*(S^2)\otimes H^*(M).
$$
For short we will say in this situation that $\pi$ is {\it c-split}.\footnote{
In some literature (see for example Thomas~\cite{Th}) this condition is called 
T.N.C.Z. (totally noncohomologous to zero), because it is equivalent 
to requiring that the inclusion of the fiber $M$ into $P$
 induce an injection on rational homology.  The paper~\cite{Mc} also discusses
situations in which the ring structure of $H^*(P)$ splits.
}
This is a deep result,
that requires the use of Gromov--Witten invariants for its proof. The results
of the present paper provide some answers to the following question:

\SmS

{\it Does any Hamiltonian fiber bundle
over a compact CW-complex  c-split? }

\SmS

A special case is when the structural group of $P\to B$ can be reduced to
a compact Lie subgroup $G$ of $\Ham(M)$.  Here c-splitting over
any base  follows from the work of
Atiyah--Bott~\cite{AB} or ours.   In this context,
one usually discusses the universal Hamiltonian $G$-bundle with fiber $M$
$$
M \;\longrightarrow\; M_G = EG\times_G M \;\longrightarrow\; BG.
$$
The cohomology of $P = M_G$ is known as the equivariant cohomology
$H_G^*(M)$ of $M$.  Atiyah--Bott show that if $G$ is a torus $T$ 
that acts in a Hamiltonian way on $M$ then the bundle $M_T \to BT$ 
is c-split.  We prove a generalization of this in
Corollary~\ref{cor:bt}.  The result for a general 
compact Lie group $G$ follows by standard arguments: 
see Corollary~\ref{cor:funct2}.   

The following theorem describes conditions on the base $B$ that
imply c-splitting.

\begin{thm}\label{thm:split}  Let $(M,\om)$ be a closed symplectic manifold,
and $M \hookrightarrow P \to B$ a bundle with structure
group $ \Ham(M)$ and 
with base a compact CW-complex $B$. Then the
rational cohomology of $P$ splits in each of the following cases:
\begin{itemize}
\item[(i)]  the base has the homotopy type of a coadjoint orbit or of a
product of spheres with at most three of
dimension~$1$;
\item[(ii)]  the base has the homotopy type of a complex blow up of a
product of complex projective spaces;
\item[(iii)] $\dim(B) \le 3$.
\end{itemize}
\end{thm}

Case $(ii)$ is a generalization of the foundational example $B = S^2$
and is proved by similar analytic methods.
The idea is to show that the map $\io: H_*(M)\to H_*(P)$ is injective by 
showing that the image $\io(a)$ in $P$ of any class  $a\in H_*(M)$
can be detected by a nonzero Gromov--Witten invariant 
of the form $n_P(\io(a), c_1,\dots,c_n; \si)$, where $c_i\in H_*(P)$ and
$\si\in H_2(P)$ is a spherical class  with nonzero image in $H_2(B)$.
The proof  should generalize to the case when all one assumes about the base
is that there is a nonzero invariant of the form $n_B( pt, pt, c_1, \ldots, c_k;
A)$: see~\cite{L} and the discussion in~\S\ref{ss:gw} below.

The proofs of parts $(i)$ and $(iii)$ 
start from the fact of c-splitting over $S^2$
and proceed  using purely topological
methods.    The following
fact about compositions of Hamiltonian bundles
is  especially   useful.  Let $M \hookrightarrow P \to
B$ be a Hamiltonian  bundle over a  simply connected base $B$ and assume that
all  Hamiltonian bundles 
over $M$ as well as over $B$ c-split. Then any Hamiltonian bundle over $P$ c-splits too.
(This fact is based on the characterization of Hamiltonian bundles in terms
of closed extensions of the symplectic form). This provides a powerful 
recursive argument which allows one to establish c-splitting over $\CP^n$
by induction on $n$,  and is an essential tool in all our 
arguments.\footnote
{
A similar property has been exploited in the context of the Halperin conjecture
discussed below:
see for example~Markl~\cite{Ma}.}

The question whether all Hamiltonian bundles over  
symplectic $4$-manifolds c-split
is still unresolved, despite 
our previous claims (cf. McDuff~\cite{BER} for example). However, even
when the base has no symplectic structure and is only a $4$-dimensional
CW-complex, our methods  still  yield some results 
about c-splitting when additional restrictions are placed on the fiber:
see \cite{LM}.
 
It is not
 clear whether one should expect that c-splitting always occurs. This question
is closely related to Halperin's conjecture, a slightly simplified version 
of which proposes that
a  fibration in the rational homotopy category whose fiber 
and base are simply connected 
 c-splits if the fiber $F$ is elliptic (that is $\pi_*(F)\otimes \Q$ has 
 finite dimension) and its rational cohomology $H^i(F)$ vanishes for odd 
 $i$.  These hypotheses imply in particular that the fiber is formal.
Clearly, the validity of Halperin's 
conjecture  with respect to a given fiber $F = (M, \om)$ 
implies that all Hamiltonian 
fibrations with that fiber are c-split.  However, note that his hypotheses 
are somewhat different from ours since  many symplectic manifolds are 
neither elliptic nor formal.
 Meier shows in~\cite{Me}~Lemma~2.5 that if $F$ is
a simply connected and
formal space such that all homotopy  
 fibrations over spheres with  fiber $F$ are 
 c-split then {\it all} fibrations with fiber $F$ and simply connected base 
are  c-split.   Since there may be homotopy fibrations that are not 
Hamiltonian,  the fact that  Hamiltonian fibrations c-split over spheres
is not enough to imply that all Hamiltonian fibrations with
simply connected and formal fiber
are c-split.   Nevertheless, Meier's result is an interesting complement 
to ours.

Although to our knowledge Halperin's conjecture is still not resolved,
there has been quite a bit of work that establishes its validity
when the fiber satisfies additional properties.    In particular,
it  holds 
when the cohomology ring $H^*(F, \Q)$ has at most $3$ even dimensional
generators (see Lupton~\cite{Lu}) or when its generators all have the same 
even dimension (see Belegradek--Kapovitch~\cite{BK}.)
In this paper we have concentrated on establishing results 
on c-splitting that hold for
all fibers $(M, \om)$.  However, there are  some simple arguments that 
apply for special $M$.  For example, 
in  section~\ref{ss:fex} we present an argument due to 
Blanchard that establishes c-splitting
when the cohomology of the fiber satisfies the hard Lefschetz
condition.  A modification due to Kedra shows that c-splitting 
holds whenever $M$ has dimension $4$.   Moreover,
the Belegradek--Kapovitch theorem has 
a Hamiltonian analog: we show in Lemma~\ref{le:dim2}
that c-splitting occurs whenever $H^*(M)$ is generated by $H^2(M)$.

In view of this, it is natural to wonder whether
 c-splitting is a purely
homotopy-theoretic property.
A c-symplectic manifold $(M, a_M)$ is defined to be a
$2n$-manifold together with a class $a_M\in H^2(M)$ 
such that $a_M^n> 0$.\footnote
{Caution: the letter ``c'' here also stands for ``cohomologically'' but
the meaning here is somewhat different from its use in the word ``c-split".
}
In view of Theorem~\ref{thm:hamchar} one could define a c-Hamiltonian
bundle  over a simply connected base manifold $B$ to be a bundle $ P\to B$
with c-symplectic fiber $(M, a_M)$ in which the symplectic class $a_M$ extends
to a class $a\in H^2(P)$.   In~\cite{All}, Allday discusses a variety of results
about symplectic torus actions, some of which do extend to the c-symplectic
case and some of which do not.  
The next lemma shows that c-splitting in general
is a geometric rather than a homotopy-theoretic property.  Its proof may be
found in \S\ref{ss:fex}.

\begin{lemma} \label{le:cham}
There is a  c-Hamiltonian bundle over $S^2$  that is not c-split.
\end{lemma}

It is also worth noting that it is essential to restrict to finite 
dimensional spaces:  c-splitting does not always hold for
``Hamiltonian" bundles with infinite dimensional fiber. (See
the footnote to Lemma~\ref{le:ev2}.)

\subsection{The homological action of $\Ham(M)$ on $M$}\label{ss:act0}

The action  $\Ham(M)\times
M \to M$  gives rise to maps  
$$
H_k(\Ham(M))\otimes H_*(M)\longrightarrow H_{k+*}(M):\quad (\phi, Z)\mapsto
\tr_\phi(Z), $$
and dually
$$
\tr^*_\phi:H_k(\Ham(M))\to \Hom(H^*(M), H^{*-k}(M)),\quad k\ge 0.
$$
In  this language, the flux of a loop
$\phi\in\pi_1(\Ham(M))$ is precisely
the element $\tr^*_\phi([\om])\in H^1(M)$. 
(Here we should use real rather than rational coefficients so that
$[\om]\in H^*(M).$)
The following result is a consequence of Theorem~\ref{thm:split}.

\begin{thm}\label{thm:act}  The maps $\tr_\phi$ and  $\tr^*_\phi$ are zero 
for all $\phi\in H_k(\Ham(M)), k> 0$.
\end{thm}

The argument goes as follows.  Recall that the cohomology ring  
of $\Ham(M)$
 is generated by elements dual to its homotopy. 
It therefore
suffices to consider the restriction of $\tr_k$ to the spherical elements $\phi$.
But in this case it is not hard to see that  the $\tr_k$
 are precisely the connecting
homomorphisms in the Wang sequence of the bundle $P_\phi\to S^{k+1}$ with clutching
function $\phi$.  These vanish because all Hamiltonian bundles 
over spheres are c-split
by part $(i)$ of Theorem~\ref{thm:split}.
Details may be found in \S\ref{ss:act}.

In particular, looking at the action on $H_0(M)$, we see that
 the point  evaluation map
$$
\ev:\Ham(M)\to M:\quad \psi\mapsto \psi(x)
$$
induces the trivial map on rational (co)homology.  
It also induces the trivial map on $\pi_1$.\footnote
{
This is a consequence of the proof of the Arnold conjecture:
see~\cite{LMP2}~\S1.3.  It is equivalent to the existence of a section of every
Hamiltonian bundle over $S^2$ and so also follows from the results
in~\cite{LMP,Mc}.}
However,
the map on $\pi_k, k> 1,$ need not be trivial. 
To see this, consider the action of $\Ham(M)$ on the symplectic frame bundle
$\SFr(M)$ of $M$ and the corresponding point evaluation maps.  
The obvious action of $\SO(3)\simeq \Ham(S^2)$ on
$\SFr(S^2) \simeq \RP^3$ induces an isomorphism 
$$
H_3(\SO(3))\cong
H_3(\SFr(S^2)),
$$
showing that these evaluation maps are not homologically trivial.  Moreover,
 its composite with the projection  $\SFr(S^2)\to S^2$ gives rise to a
nonzero map 
$$
\pi_3(SO(3)) = \pi_3(\Ham(S^2))\to \pi_3(S^2).
$$
Thus the corresponding Hamiltonian fibration over $S^4$ with fiber $S^2$, though
c-split, does not have a section. 

 Note, however, that  the evaluation map
$$
\pi_{2\ell}(X^X) \stackrel{\ev}\to \pi_{2\ell}(X) \to H_{2\ell}(X,\Q),\quad \ell > 0,
$$
is  always zero,
if $X$ is a finite CW complex and $X^X$ is its space
of self-maps.
Indeed, because the cohomology ring $H^*(X^X,\Q)$ is freely generated
by elements dual to $\pi_*(X^X)\otimes \Q$,  there would otherwise be an element
$a\in H^{2\ell}(X)$  that would pull back to an element of infinite order in the
cohomology ring of the $H$-space $X^X$. Hence $a$ itself would have to have
infinite order, which is impossible.
  A more delicate argument shows that
the integral evaluation $\pi_{2\ell}(X^X)\to H_{2\ell}(X,\Z)$ is zero: 
see~\cite{Go}. 
\MS

  By Lemma~\ref{le:auttriv}, a Hamiltonian automorphism  of the product
Hamiltonian bundle  $B\times M \to B$ is simply a map $B \to B\times
\Ham(M)$ of the form $b\mapsto (b, \phi_b)$.   If $B$ is a closed manifold
we will see that Theorem~\ref{thm:act} implies
that any Hamiltonian automorphism  of the product  bundle acts as the identity
map on the rational cohomology of $B\times M$: see
Proposition~\ref{prop:equiv}.   The natural generalization
of this result would claim that a Hamiltonian automorphism of a 
bundle $P$ acts
as the identity map on the rational cohomology of $P$.  We do not know yet
whether this is true in  general. However, we can show that it is closely
related to the c-splitting of Hamiltonian bundles. Thus we can establish it only
under conditions similar  to the conditions under which c-splitting holds. See 
Proposition~\ref{prop:automorphisms} below.

\subsection{Implications for general symplectic bundles}\label{ss:gen}

Consider the Wang sequence for a
symplectic bundle $\pi:P\to S^2$ with clutching map $\phi\in
\pi_1(\Symp(M))$: 
$$
\cdots \; \longrightarrow  H^k(M)\stackrel{\p}{\longrightarrow} H^{k-1}(M)
\stackrel{u}{\longrightarrow} H^{k+1}(P)
\stackrel{restr}{\longrightarrow} H^{k+1}(M)
\longrightarrow \; \cdots
$$
Here the map $u$ may be realized in de Rham cohomology by choosing any
extension of a given closed form $\al$ on $M$ and then wedging it with the
pullback of a normalized area form on the base.  Further, as pointed out above,
the boundary map $\p = \p_\phi$ is just $\tr^*_\phi$.   Thus the bundle is 
Hamiltonian if and only if $\tr^*_\phi([\om]) = \p([\om])= 0$.  In
the Hamiltonian case Theorem~\ref{thm:split} implies that $\p$ is
identically $0$.  In the general case, we know that the map $\p: H^*(M)\to
H^{*-1}(M)$ is a derivation: i.e.  
$$
\p(ab) = \p(a)b + (-1)^{deg(a)} a \p(b).
$$
The following result is an easy consequence of the fact that the action of
$\pi_2(\Ham(M)) = \pi_2(\Symp(M))$ on $H^*(M)$ is trivial.

\begin{prop}\label{prop:p}  The boundary map $\p$ in the Wang
rational cohomology sequence of a symplectic bundle over $S^2$ has the
property that $\p\circ\p = 0$. \end{prop}

The proof is given in \S\ref{ss:nonH}.
The above result holds trivially when $\phi$ corresponds to a smooth (not
necessarily symplectic)  $S^1$-action since then $\p$ is given in deRham
cohomology by contraction $\iota_X$ by the generating vector field $X$.  
Moreover, the authors know of no smooth bundle over $S^2$ for which the
above proposition does not hold, though it is likely that they exist.  By the remarks
in~\S\ref{ss:nonH} such a bundle would have no extension over $\CP^2$.

One consequence is the following result about the
boundary map $\p = \p_\phi$  in the case when the loop $\phi$ is far from
being Hamiltonian.  Recall (e.g. from~\cite{LMP}) that $\pi_1(\Ham(M))$ is
included in (but not necessarily equal to) the kernel of the evaluation map 
$\pi_1(\Symp(M))\to \pi_1(M)$.  Any loop whose evaluation is homologically
essential can therefore be thought of as ``very nonHamiltonian".

\begin{cor}\label{cor:vnH}  $\kerr\, \p = \im\, \p$ if and only if
the image of $\phi$ under the evaluation map $\pi_1(\Symp(M))\to
H_1(M,\Q)$ is nonzero. 
\end{cor}

A similar result was obtained by Allday concerning
$S^1$ actions on c-symplectic manifolds: see statement (d) in~\cite{All}.
He was considering manifolds $M$ that
satisfy the weak Lefschetz condition, i.e. manifolds such that
$$
\wedge [\om]^{n-1}: \quad H^1(M,\R) \to
H^{2n-1}(M,\R)
$$
 is an isomorphism, in which case every nonHamiltonian loop is 
``very nonHamiltonian.''

\section{The characterization of Hamiltonian bundles} \label{se:hamchar}

This section contains proofs of the basic
results on the existence and classification of
Hamiltonian bundles, namely Theorems~\ref{thm:hamchar}
and~\ref{thm:hamchar3} and
Proposition~\ref{prop:hamchar2}.

\subsection{Existence of Hamiltonian structures}\label{ss:hamchar}

We begin by giving a proof of
Theorem~\ref{thm:hamchar} using as little analysis as possible.
We will repeat some of the arguments in~\cite{MS} Ch. 6 for the sake of clarity.
The main new point is the replacement of the Guillemin--Lerman--Sternberg
construction of the coupling form by the more topological 
Lemma~\ref{le:ev2}.

Geometric proofs (such as those in~\cite{MS}) apply  when $P$ and $B$
are smooth manifolds and  $\pi$ is a smooth
surjection.  However, as the following lemma makes clear, this 
is no restriction.

\begin{lemma}\label{le:smooth}  Suppose that $\pi: Q\to W$ is a locally
trivial  bundle over a finite
CW complex $W$ with compact fiber $(M, \om)$ and suppose that 
the structural group $G$
is equal either to $\Symp (M,\om)$   or to $\Ham(M,\om)$.  Then 
there is a smooth bundle $\pi:P\to B$ as above with structural group $G$ 
and a homeomorphism $f$ of $W$ onto a closed subset of $B$ such that
$\pi:Q\to W$ is homeomorphic to the pullback bundle $f^*(P)\to W$.
\end{lemma}
\proof{}   Embed
$W$ into  some Euclidean space and let $B$ be a suitable small neighborhood of
$W$.  Then $W$ is a retract of $B$ so that the classifying map $W\to
B\,G$ extends to $B$.  It remains to approximate this map $B \to B\,G$ by a
smooth map. \QED

First let us  sketch the proof of
Theorem~\ref{thm:hamchar} when the base is simply connected.
We use the minimum amount of geometry: nevertheless to get a relation
between the existence of the class $a$ and the 
structural group it seems necessary 
to use the idea of a symplectic connection.   We begin with an easy lemma.

\begin{lemma}\label{le:conn}  Let $P\to B$ be a symplectic bundle with
closed  connection form $\tau$.  Then the holonomy of the corresponding
connection $\Nabla_\tau$ round any contractible loop in $B$ is Hamiltonian.
 \end{lemma} 
\proof{} It suffices to consider the case when $B = D^2
$.  Then the bundle
$\pi:P\to D^2$ is symplectically trivial and 
so may be identified with the product
$D^2\times M$ in such a way that the symplectic 
form on each fiber is simply
$\om$.  
Use this trivialization to identify the 
 holonomy round the loop
$s\mapsto e^{2\pi i s}\in \p D^2$  with a family of
diffeomorphisms $\Phi_s:M\to M, s\in [0,1]$. 
Since this holonomy is simply the flow along 
the null directions (or characteristics) of the closed
form $\tau$ on the hypersurface $\p P$, 
a standard calculation shows that the
$\Phi_s$ are symplectomorphisms.
 Given a $1$-cycle $\de:S^1\to
M$ in the fiber $M$ over $1\in \p D^2$, 
consider the closed $2$-cycle that is the
union of the following two cylinders:  \begin{eqnarray*} 
c_1: [0,1]\times S^1\to \p D^2\times M: & &
(s,t)\mapsto (e^{2\pi i s},\Phi_s(\de(t))),\\
c_2: [0,1]\times S^1\to  1\times M: & &
(s,t)\mapsto (1,\Phi_{1-s}(\de(t))).
\end{eqnarray*}
This cycle is obviously contractible.  Hence, 
$$
\tau(c_1) = -\tau(c_2) = \Flux(\{\Phi_s\})(\de).
$$
But $\tau(c_1) = 0$ since the characteristics of $\tau|_{\p P}$ are tangent to
$c_1$.  Applying this to all $\de$, we see that  the holonomy round $\p D^2$ has
zero flux and so is Hamiltonian.\QED

\begin{lemma}\label{le:shamchar}  If $\pi_1(B) = 0$ then a symplectic bundle
$\pi:P\to B$
 is Hamiltonian if and only if the class $[\om_b]\in H^2(M)$ extends to $a\in
H^*(P)$.  \end{lemma}
\proof{}  Suppose first that the class $a$ exists.  By Lemma~\ref{le:smooth},
we can work in the smooth category.   Then Thurston's convexity argument
allows us to construct a closed connection form $\tau$ on $P$ and hence a
 horizontal distribution  $Hor_\tau$.  The previous lemma shows that the
holonomy  around every contractible loop in $B$  is Hamiltonian.  Since $B$ is
simply connected, the holonomy round {\it all} loops is Hamiltonian.  Using this,
it is easy to reduce the structural group of  $P\to B$  to $\Ham(M)$.  For more
details, see~\cite{MS}.

Next, suppose that the bundle is Hamiltonian.  We need to show that the fiber
symplectic class extends to $P$.  The proof in~\cite{MS} does this by
the method of Guillemin--Lerman--Sternberg and constructs  
a closed connection form $\tau$ (called the
coupling form) starting from a  connection with Hamiltonian holonomy.
This construction uses the curvature of the connection and is quite analytic.  
In contrast, we shall now use topological arguments to reduce to the cases
$B = S^2$ and $B= S^3$.  These cases are then dealt with by  elementary
arguments.

Consider the Leray--Serre cohomology spectral sequence for $M\to P\to B$.
Its $E_2$ term is a product:  $E_2^{p,q} = H^p(B)\otimes H^q(M)$.  (Here 
$H^*$ denotes cohomology over $\R$.)  We need to show that the class
$[\om] \in E_2^{0,2}$ survives into the $E_\infty$ term, 
which happens if and only if
it is in the kernel of
the two differentials  $d_2^{0,2}, d_3^{0,2}$.  
Now
$$
d_2^{0,2}: H^2(M) \to H^2(B)\otimes H^1(M)
$$
is essentially the same as the flux homomorphism.  More precisely, 
if $c:S^2\to B$ represents some element (also called $c$) in $H_2(B)$, then the
pullback of the bundle $\pi: P\to B$  by $c$ is a bundle over $S^2$ that is
determined by a loop $\phi_c\in \pi_1(\Ham(M))$ that is well defined up to
conjugacy.  Moreover, for each $\la\in H_1(M)$,
$$
d_2^{0,2} ([\om]) (c, \la) = \tr^*_{\phi_c}(\la),
$$
where $\tr^*$ is as in \S\ref{ss:act0}.
Hence $d_2^{0,2}([\om]) = 0$ because $\phi_c$ is Hamiltonian.  

To deal with $d_3$ observe first that because
 the inclusion of the
$3$-skeleton $B_3$  into $B$ induces an injection 
$H^q(B)\to H^q(B_3)$ for $q\le 3$,
$d_3^{0,2}([\om])$ vanishes in the spectral sequence for $P\to B$
if  it vanishes for the pullback bundle over $B_3$.
Therefore we may
suppose that $B$ is a $3$-dimensional CW-complex whose $2$-skeleton
$B_2$ is a  wedge of $2$ spheres.  (Recall that $\pi_1(B)=0$.)
Further,  we can choose the cell decomposition so that the first $k$ $3$-cells
span the kernel of the boundary map $C_3 \to C_2$ in the cellular chain
complex of $B_3$.  Because $H_2(B_2) = \pi_2(B_2)$, the attaching maps of
these first $k$-cells are null homotopic.  Hence there is  a wedge $B'$ of
$2$-spheres and $3$-spheres and a map $B'\to B_3$ that induces a
surjection on $H_3$. It therefore suffices to show that 
$d_3^{0,2}([\om])$ vanishes in the pullback bundle over $B'$.
This will clearly be the case if it vanishes 
in every Hamiltonian bundle over
$S^3$.

Now, a Hamiltonian fiber bundle over $S^3$ is determined by a map 
$$
I^2/ \p I^2= S^2 \to
\Ham(M):\quad (s,t)\mapsto \phi_{s,t},
$$ 
and it is easy to see that $d_3^{0,2}([\om])= 0$ exactly when the
 the evaluation map
$$
\ev_x:\Ham(M)\to M: \phi\mapsto \phi(x)
$$
takes $\pi_2(\Ham(M))$ into the kernel of $\om$.

The result now follows from Lemma~\ref{le:ev2} below.\QED

\begin{lemma}\label{le:ev2}  Given a smooth map
$
\Psi: (I^2, \p I^2) \to
(\Ham(M), id)
$
and $x\in M$, let $\Psi^x :(I^2, \p I^2)\to M$ be the composite of $\Psi$ with
evaluation at $x$.   Then
$$
\int_{I^2}(\Psi^x)^*\om = 0,
$$
for all $x\in M$.
\end{lemma}
\proof{} For each $s,t$ let $X_{s,t}$ (resp. $Y_{s, t}$) be the Hamiltonian vector
field on $M$ that is tangent to the flow of the isotopy $s\mapsto \Psi^x(s, t)$,
(resp. $t\mapsto \Psi^x(s, t)$.)
Then
$$
 \int_{I^2} (\Psi^x)^*\om = \int\!\!\int \om(X_{s,t}(\Psi^x(s,t)),
Y_{s,t}(\Psi^x(s,t)))\, ds dt. $$
The first observation is that this integral is a 
constant $c$ that is independent of  $x$, since
the  maps $\Psi^x:S^2\to M$ are all homotopic.  Secondly, 
recall that for any Hamiltonian
vector fields $X,Y$ on $M$
$$
\int_M \om(X,Y)  \om^n = n \int_M \om(X,\cdot) \om(Y,\cdot) \om^{n-1} = 0,
$$
since $\om(X,\cdot), \om(Y,\cdot)$ are exact $1$-forms.  
Taking $X_{s,t} = X_{s,t}(\Psi^x(s,t))$ and similarly for $Y$, we have
$$
\int c\,\om^n = \int_{I^2} \left(\int_ M \om(X_{s,t},Y_{s,t})\, \om^n\right)
ds\, dt\; =\; 0. $$
Hence $c=0$.

This lemma can also be proved by purely topological methods.  In fact, 
as remarked in the discussion after Theorem~\ref{thm:act}, the
evaluation map $\pi_2(X^X)\to H_2(X,\R)$ vanishes for any finite CW
complex $X$.\footnote
{
The following example due to Gotay {\it et al}~\cite{GLSW}
demonstrates the importance of this finiteness hypothesis.
 Let $\Hh$ be a complex Hilbert space
with unitary group $\U(\Hh)$ and consider the exact sequence 
$
S^1\to \U(\Hh) \to \PU(\Hh),
$ 
where $\PU(\Hh)$ is the projective unitary group. 
 Since $\U(\Hh)$ is contractible, $\PU(\Hh)\simeq
\CP^\infty$.  Since $\PU(\Hh)$ can be considered as a subgroup of
the symplectomorphism group of
$\CP(\Hh)$, the generator $\phi$ of $\pi_2(\PU(\Hh))$ gives
rise to a ``symplectic" fibration
$
\CP(\Hh) \to P_\phi \to S^3,
$
which is ``Hamiltonian'' because $\CP(\Hh)\simeq \CP^\infty$ is simply connected.  
It is easy to check that the evaluation map
$
\PU(\Hh) \to \CP(\Hh):  A\mapsto A(x)
$
is a homotopy equivalence.  Hence Lemma~\ref{le:shamchar} fails in this case.}

\QED

This completes the proof of Theorem~\ref{thm:hamchar} 
for simply connected bases.  

 \begin{lemma} Theorem~\ref{thm:hamchar}
 holds  for all $B$.
\end{lemma}
\proof{} Suppose that $\pi_B: P\to B$ is Hamiltonian.  It is classified by a map 
$B\to B\,\Ham(M)$.  Because $B\,\Ham(M)$ is simply connected 
this factors through a map $C\to B\,\Ham(M)$,
where $C$ is obtained by collapsing 
the $1$-skeleton of $B$ to a point.  In particular condition $(i)$ is satisfied. 
To verify $(ii)$, let $\pi_C:
Q\to C$ be
 the corresponding Hamiltonian bundle, so that there
 is a commutative diagram
$$
\begin{array}{rcc} P &\to & Q\\
\pi_B\downarrow & &\pi_C \downarrow\\
B &\to & C = B/B_1.
\end{array}
$$
 Lemma~\ref{le:shamchar} applied  to $\pi_C$ 
tells us that there is a class $a_C\in
H^2(Q)$ that restricts to $[\om]$ on the fibers.  Its pullback to $P$ is the desired
class $a$.

Conversely, suppose that conditions $(i)$ and $(ii)$ are satisfied.  
By $(i)$, the classifying map
$B\to B\,\Symp(M)$ factors through a map $f: C\to B\,\Symp(M)$,
where $C$ is as above. This map $f$ depends on the choice of 
a symplectic trivialization of
$\pi$ over the $1$-skeleton $B_1$ of $B$.  We now show that $f$ can be chosen
so that $(ii)$  holds for the associated symplectic bundle $Q_f\to C$.

As in the proof of Lemma~\ref{le:shamchar}. we need to show that the
differentials  $(d_C)_2^{0,2}, (d_C)_3^{0,2}$  in the spectral sequence for $Q_f\to
C$ both vanish on $[\om]$.
  Let 
$$
\cdots \to C_k(B)\stackrel{\p}\to C_{k-1}(B)\to\cdots
$$
be the cellular chain complex for $B$, and
choose $2$-cells $e_1,\dots, e_k$
in $B$ whose attaching maps $\al_1, \dots, \al_k$ form a basis over $\Q$ for the
image of $\p$ in $C_1(B)$. Then the obvious maps $C_k(B) \to C_k(C)$ (which are
the identity for $k > 1$) give rise to an isomorphism 
$$
 H_2(B,\Q) \,\bigoplus\, \oplus_i\Q[ e_i] \cong H_2(C, \Q).
$$
By the naturality of  spectral sequences, the vanishing of 
$(d_B)_2^{0,2}([\om])$ implies that
$
(d_C)_2^{0,2}([\om])$ vanishes on all cycles in $H_2(C, \Q)$ coming from $H_2(B,
\Q)$.  Therefore we just need to check that it vanishes on the cycles $e_i$.
For this, we have to choose the trivialization over $B_1$ so that its pullback by
each $\al_i$ gives rise to a  Hamiltonian bundle over $e_i$.  
For this it would
suffice  that its pullback by each $\al_i$ is the ``natural
trivialization", i.e the one that extends over the $2$-cell $e_i$. 
To arrange this, choose any symplectic
trivialization over $B_1 = \vee_j \ga_j$. Then comparing this with the natural
trivialization gives rise to a homomorphism
$$
\Phi: \oplus_i \Z e_i \to \pi_1\Symp(M,\Z) \stackrel{\Flux}\longrightarrow
H^1(M,\R). $$
 Since the boundary map 
$\oplus_i \Z e_i\to C_1(B)\otimes
\Q$ is injective, we can now change the chosen trivializations over the $1$-cells
$\ga_j$ in $B_1$  to make $\Phi = 0$.

This ensures that $d_2^{0,2} = 0$ in the bundle over $C$.  
Since the map $H^q(C)\to H^q(B)$ is an isomorphism when $q\ge 3$, the vanishing
of $d_3^{0,2}$ for $B$ implies that it vanishes for $C$.  Therefore $(ii)$ holds for
$Q\to C$.  By the previous result, this implies that the structural group of $Q\to C$ 
reduces to $\Ham(M)$.  Therefore, the same holds for $P\to B$.\QED

In the course of the above proof we established the following useful result.

\begin{cor}\label{cor:hamchar2}  Let $C$ be the CW complex obtained by collapsing the
$1$-skeleton of $B$ to a point and $f:B\to C$ be the obvious map.  Then any Hamiltonian
bundle $P\to B$ is the pullback by $f$ of some Hamiltonian bundle over $C$.
\end{cor}

Theorem~\ref{thm:hamchar} shows that there are two obstructions to the
existence of a Hamiltonian structure on a symplectic bundle.  Firstly, the
bundle must be symplectically trivial over the $1$-skeleton $B_1$, and
secondly the symplectic class on the fiber must extend.  The first obstruction
obviously depends on the $1$-skeleton $B_1$ while the second, in principle,
depends on its $3$-skeleton (since we need $d_2$ and $d_3$ to vanish on
$[\om]$).  However, in fact, it only depends on the $2$-skeleton, as is shown in
the next lemma.

\begin{lemma}\label{le:ham2c}  Every symplectic bundle 
over a $2$-connected
base $B$ is Hamiltonian.
\end{lemma}
\proof{}  We give two proofs.  First, observe that
as in Lemma~\ref{le:shamchar} 
we just have to show that $d_3^{0,2}([\om]) = 0$.  
The arguments of that lemma apply to show that this is the case.

Alternatively, let $\TS$ (resp. $\TIH$) denote the universal cover
of the group $\Symp_0 = \Symp_0(M,\om)$ (resp $\Ham(M)$),  and  set $\pi_S =
\pi_1(\Symp_0)$ so that there are fibrations $$
\TIH\to \TS \stackrel{\Flux}\to H^1(M,\R),\quad   B(\pi_S)\to B\,\TS \to
B\,\Symp_0. $$
The existence of  the first fibration shows that  $\TIH$  is homotopy equivalent to $\TS$ so
that $B\,\TIH\simeq B\,\TS$, while the second implies that there is a fibration
$$
B\,\TS \to B\,\Symp_0\to K(\pi_S,2),
$$
where $K(\pi_S,2)$ is an Eilenberg--MacLane space.  A symplectic
bundle over $B$ is equivalent to a homotopy class of maps $B\to B\,\Symp_0$. If
$B$ is $2$-connected,  the composite $B\to B\,\Symp_0 \to K(\pi_S,2)$ is null
homotopic, so that the map $B\to B\,\Symp_0$ lifts to $B\,\TS$ and hence to the
homotopic space $B\,\TIH$. Composing this map $B\to B\,\TIH$ with the projection
$B\,\TIH \to B\,\Ham$ we get a 
Hamiltonian structure on the given bundle over $B$.  

Equivalently, use the existence of the fibration 
$\TIH\to \TS \to H^1(M,\R)$ to deduce that the subgroup $\pi_1(\Ham)$
of $\TIH$ injects into $\pi_1(\Symp_0)$.  This implies that the relative
homotopy groups $\pi_i(\Symp_0, \Ham)$ vanish for $i > 1$, so that
$$
\pi_i(B\,\Symp_0, B\,\Ham) = \pi_{i-1}(\Symp_0, \Ham) = 0,\quad i>2.
$$
The desired conclusion now follows by obstruction theory.
 \QED

The second proof does not  directly use the 
sequence $0 \to \Ham \to \Symp_0 \to H^1/\Ga_\om\to 0$ since the flux group
$\Ga_\om$ may not be a discrete subgroup of $H^1$.

\subsection{The classification of Hamiltonian structures}\label{ss:hamstr}

The previous subsection discussed the question of the existence of
Hamiltonian structures on a given bundle.  We now look at the problem of
describing and classifying them.  We begin by proving 
Lemma~\ref{le:normalized} that states that any 
closed extension of the fiber class can be normalized.
\MS

\NI
{\bf Proof of Lemma~\ref{le:normalized}}\,\,  
Let $\pi:P\to B$ be a symplectic bundle satisfying the
conditions of Theorem~\ref{thm:hamchar}
and fix an identification of $(M, \om )$ with $(M_{b_0}, \om_{b_0})$. 
 Let $a$ be any closed extension
of $[\om]$, $\ga_1, \ldots, \ga_k$ be a set of generators of the first rational
homology group of $B$, $\{c_i\}$ the dual basis of $H^1(B)$ and $T_1, \ldots,
T_k$  symplectic trivializations round the $\ga_i$. Assume for the moment
that each class $f(T_i,a) \in H^1(M_{b_0}) = H^1(M)$ has an extension
$\tilde{f}(T_i,a)$ to $P$. Subtracting from
$a$ the class $\sum_{i=1}^{k} \pi^*(c_i) \cup \tilde{f}(T_i,a)$, we get a 
closed extension $a'$ whose corresponding classes  $f(T_i,a')$ belong
to $\Ga_{\om}$.

   There remains to prove that the extensions of the $f(T_i,a)$'s exist 
in Hamiltonian bundles. It is enough to prove that the fiber inclusion
$M \to P$ induces an injection on the first homology group.
One only needs to
prove this over the $2$\/-skeleton $B_2$ of $B$, and, by
Corollary~\ref{cor:hamchar2} we can assume as well that $B_2$ is a wedge of 
$2$\/-spheres. Hence this is a consequence of the easy fact that the evaluation of a
Hamiltonian loop  on a point of $M$ gives a $1$-cycle of $M$ that is trivial in rational
homology, i.e. that the differential $d_2^{0,1}$ vanishes in the cohomology spectral
sequence for $P\to B$:  see for instance \cite{LMP} where this is proved by
elementary methods.  \QED

The next result extends Lemma~\ref{le:conn}. 

\begin{lemma}\label{le:conn2}  Let $P\to B$ be a symplectic bundle with
 a given symplectic trivialization of $P$ over $B_1$, and let $a\in H^2(P)$
be a normalized extension of the fiber symplectic class.   Then the restriction
of $a$ to $\pi^{-1}(B_1)$ defines and is defined by a
homomorphism $\Phi$ from $\pi_1(B)$ to $\Ga_{\om}$.
 \end{lemma} 
\proof{} As in Lemma~\ref{le:conn}, we can use the given trivialization to
identify the  holonomy round some loop
$s\mapsto \ga(s)\in B_1$  with a family of
symplectomorphisms $\Phi_s^\ga:M\to M, s\in [0,1]$. 
 Given a $1$-cycle $\de:S^1\to
M$ in the fiber $M$ over $1\in \p D^2$, 
consider the closed $2$-cycle $C(\ga,\de) = c_1\cup c_2$ as before.
Since $\tau(c_1) = 0$,
$$
\tau(C(\ga,\de)) = \tau(c_2) = -\Flux(\{\Phi_s^\ga\})(\de).
$$
If we now set
$$
\Phi(\ga) = - \Flux(\{\Phi_s^\ga\}),
$$
it is  easy to check that $\Phi$
is a homomorphism. Its values are in $\Ga_{\om}$ by the definition of
normalized extension classes.  The result follows. \QED

The next task is to prove Theorem~\ref{thm:hamchar3} that characterizes
Hamiltonian structures.  Thus we need to understand the homotopy classes
of lifts $\tg$ of the classifying map $g: B \to B\,\Symp(M, \om)$ of the
underlying symplectic bundle to $B\,\Ham(M)$.
We first consider the intermediate lift  $\hg$ of $g$ into
 $B\,\Symp_0(M,\om)$.  
In view of the fibration sequence
$$
\pi_0(\Symp) \to B\,\Symp_0 \to B\,\Symp \to B(\pi_0(\Symp))
$$
in which each space is mapped to the homotopy fiber of the subsequent map,
a map $g: B\to B\,\Symp$ lifts to $\hg: B\to B\,\Symp_0$ if and only if the
symplectic bundle given by $g$ can be trivialized over the
$1$-skeleton $B_1$ of $B$.  Moreover such lifts are in bijective correspondence
with the elements of $\pi_0(\Symp)$ and so correspond to an identification (up to
symplectic isotopy) of $(M, \om)$ with the fiber $(M_{b_0}, \om_{b_0})$ at the
base point $b_0$.  (Recall that $B$ is always assumed to be connected.) 

To understand the full lift $\tg$, recall 
the exact sequence
$$
\{id\} \longrightarrow \Ham(M,\om) \longrightarrow \Symp_0(M,\om) 
\stackrel{\Flux}\longrightarrow H^1(M,\R)/\Ga_{\om}\longrightarrow
\{0\}.\qquad\quad(*) $$
If $\Ga_{\om}$ is discrete, then
the space $H^1(M,\R)/\Ga_{\om}$ is homotopy equivalent to a torus 
and we can investigate the liftings $\tg$ by homotopy theoretic arguments
about the fibration 
$$
H^1(M,\R)/\Ga_{\om} \to B\,\Ham(M,\om) \to B\,\Symp_0(M,\om).
$$
However, in general, we need to argue more geometrically.

Suppose that  a symplectic bundle $\pi:P\to B$
is given that satisfies the conditions of Theorem~\ref{thm:hamchar}.
Fix an identification of $(M, \om )$ with $(M_{b_0}, \om_{b_0})$.
We have to show that  lifts from $B\,\Symp_0$
to $B\,\Ham$ are in bijective correspondence with equivalence classes of
normalized extensions  $a$ of the fiber symplectic class. 
By Theorem~\ref{thm:hamchar} and Lemma~~\ref{le:normalized}, 
there is a lift if and only if there is
 a normalized extension class $a$.  Therefore, it remains to show that 
the equivalence relations correspond.
The essential reason why this is true is that the induced map
$$
\pi_i(\Ham(M,\om)) \to \pi_i(\Symp_0(M, \om))
$$
is an injection for $i = 1$ and an isomorphism for $i > 1$.  This, in turn, follows
from the exactness of the sequence $(*)$.

Let us spell out a few more details, first
when $B$ is simply connected.  Then  the classifying map from the
$2$-skeleton $B_2$ to $B\,\Symp_0$ has a lift to $B\,\Ham$ if and only if the
image of the induced map 
$$
\pi_2(B_2) \to \pi_2(B\,\Symp_0(M)) = \pi_1\Symp_0(M)
$$
lies in the kernel of the flux homomorphism 
$$
\Flux: \; \pi_1(\Symp_0(M)) \longrightarrow  \Ga_{\om}.
$$
Since $\pi_1(\Ham(M,\om))$ injects into $\pi_1(\Symp_0(M, \om))$
there is only one such lift up to homotopy.  Standard arguments now show that
this lift can be extended uniquely to the rest of $B$.  Hence in this case 
there is a unique lift.  Correspondingly there is a unique equivalence class of
extensions $a$.

Now let us consider the general case.
We are given a map $g: B\to B\,\Symp_0$ and want to identify the different
homotopy classes of liftings of $g$ to $B\,\Ham$.  Let $C = B/B_1$ as above. 
By Corollary~\ref{cor:hamchar2}, there is a
symplectic trivilization $T$ over $B_1$ 
that is compatible with the given identification of 
the base fiber and induces a map
$C\to B\,\Symp_0$ which lifts to $B\,\Ham$.
Since this lifting $g_{T,C}$ of $C\to B\,\Symp_0$ is unique,
each isotopy class $T$ of such trivializations over $B_1$ gives rise to a unique
homotopy class $g_T$  of maps $B\to B\,\Ham$, namely
$$
g_T:\quad B\longrightarrow  C\stackrel{g_{T,C}} \longrightarrow  B\,\Ham.
$$
Note that $g_T$ is a lifting of $f$ and that every lifting occurs this way.

Standard arguments show 
 that two such isotopy classes  differ by a homomorphism
$$
\pi_1(B)\longrightarrow  \pi_1(\Symp_0).
$$
Moreover,  the corresponding 
maps $g_T$ and $g_T'$
are homotopic if and only if $T$ and $T'$ differ by a homomorphism
with values in $\pi_1(\Ham)$.   Thus homotopy classes of liftings of $g$ to
$B\,\Ham$ are classified by homomorphisms
$\pi_1(B) \to \Ga_{\om}$.  By Lemma~\ref{le:conn2} these homomorphisms
are precisely what defines the equivalence classes of extensions $a$.\QED

\section{Properties of Hamiltonian bundles}

The key to extending results about Hamiltonian bundles over $S^2$ 
to higher dimensional bases  is  their functorial
properties, in particular their behavior under composition.  Before
discussing this, it is useful to establish the fact that this class of bundles is
 stable under  small perturbations of the symplectic form on $M$.

\subsection{Stability}\label{ss:stab}

Moser's argument implies that for every 
symplectic structure $\om$ on $M$ there
is a Serre fibration 
$$
\Symp(M,\om) \longrightarrow \Diff(M) \longrightarrow  \Ss_{\om},
$$
where $\Ss_{\om}$ is the space of all symplectic structures on $M$ that are
diffeomorphic to $\om$.
At the level of classifying spaces, this gives a homotopy fibration
$$
\Ss_{\om} \hookrightarrow  B \Symp(M, \om) \longrightarrow B \Diff(M).
$$
Any smooth fiber bundle $P\to B$ with fiber $M$ is
classified by a map $B \to B\,\Diff(M)$, and isomorphism classes of symplectic
structures on it with fiber $(M, \om)$ correspond to homotopy classes of
sections of the associated fibration $W(\om) \to B$ with fiber $\Ss_{\om}$.
We will suppose that $\pi$ is described by  a finite set of local
trivializations $T_i: \pi^{-1}(V_i) \to V_i\times M$ with
the transition functions  
$\phi_{ij}: V_i\cap V_j \to \Diff(M)$.

\begin{lemma}\label{le:sstr}
 Suppose that $M\to P\to B$ is a smooth 
fibration constructed from a cocycle $(T_, \phi_{ij})$ with the following 
property:  there is a symplectic form $\om$ on $M$ such that for each 
$x\in M$ the convex hull of the finite set of forms 
$$
\{\phi_{ij}^*(\om): x\in V_i\cap V_j\}
$$
lies in the set $\Ss_\om$ of symplectic forms diffeomorphic to $\om$.
Then $(M,\om) \to P\to B$ may be given the structure of a symplectic bundle.
\end{lemma}
\proof{}  It suffices to construct a section $\si$ of $W(\om)\to B$. 
 The hypothesis 
implies that for each $x$ the convex hull of the set of forms 
$T_i(x)^*(\om), x\in V_i,$ lies in the fiber of $W(\om)$ at $x$.
Hence we may take
$$
\si(x) = \sum_i \rho_i  T_i^*(\om),
$$
where $\rho_i$ is a partition of unity subordinate to the cover $V_i$.\QED

\begin{cor}\label{cor:stab0} 
Let $P\to B$ be a symplectic bundle with closed fiber $(M, \om)$ and 
compact base $B$.  There is an open neighborhood 
$\Uu$ of $\om$ in the space $\Ss(M)$ of all symplectic forms on $M$
such that,  for all $\om'\in U$, the structural group of $\pi:P\to B$ may be
reduced to $\Symp(M,\om')$.  
\end{cor}
\proof{} Trivialize $P\to B$ so that $\phi_{ij}^*(\om) = \om$ for all 
$i,j$.  Then the hypothesis of the lemma is satisfied for all $\om'$ 
sufficiently close to $\om$ by the openness of the symplectic 
condition.\QED

Thus the set
$\Ss_\pi(M)$ of symplectic forms on
$M$, with  respect to which $\pi$ is symplectic, is open.
The aim of this section is to show that a corresponding statement is true
for  Hamiltonian bundles.  The following example
shows that the Hamiltonian property need not survive under {\it large}
perturbations of $\om$ because condition $(i)$ in 
Theorem~\ref{thm:hamchar} can
fail.  However, it follows from the proof of stability that 
condition $(ii)$ never fails under any perturbation.

\begin{ex}\label{ex:stab} \rm  Here is an example of a smooth family of
symplectic bundles  that is Hamiltonian at all times $0 \le t < 1$ but is
nonHamiltonian at time $1$. Let $h_t, 0\le t\le 1$, be a family of
diffeomorphisms of $M$ with $h_0 = id$ and define $$
Q = M\times [0,1]\times [0,1]/(x,1,t) \equiv (h_t(x), 0,t).
$$
Thus we can think of $Q$ as  a family of  bundles $\pi: P_t\to S^1$ with monodromy
$h_t$ at time $t$.  Seidel~\cite{Sei} has shown that there are smooth families of
symplectic forms $\om_t$ and diffeomorphisms $h_t\in \Symp(M,\om_t)$ for $t
\in [0,1]$ such that $h_t$ is not in the identity component of $\Symp(M,\om_t)$
for $t=1$ but is in this component for $t < 1$.  For such $h_t$ each bundle
$P_t\to S^1$ is symplectic.  Moreover, it  is symplectically trivial and hence
Hamiltonian for $t < 1$ but is nonHamiltonian at $t=1$.  This example shows why
the verification of condition $(i)$ in the next proof is somewhat delicate.
\end{ex}

 \begin{lemma}\label{le:st2}  A Hamiltonian bundle $\pi:P\to B$ is stable
if and only if the restriction map
$H^2(P)\to H^2(M)$ is surjective.\end{lemma}

\proof{}   If $\pi:P\to B$ is
Hamiltonian with respect to $\om'$ then by  Theorem~\ref{thm:hamchar}
 $[\om']$ is in the image of  
$H^2(P)\to H^2(M)$.  If $\pi$ is stable, then $[\om']$ fills out a neighborhood of $[\om]$
which implies surjectivity.  Conversely, suppose that we have surjectivity.  Then 
the second condition of Theorem~\ref{thm:hamchar} is satisfied. 
To check $(i)$ let $\ga: S^1\to B$ be a loop in $B$ and suppose that
$\ga^*(P)$ is identified symplectically with the product bundle
$S^1\times (M,\om)$.  
Let $\om_t, 0\le t\le \eps,$ be a (short) smooth path with $\om_0 =
\om$.   Then, because $P\to B$ has the structure of an $\om_t$-symplectic bundle for each
$t$, each fiber $M_b$ has a corresponding  smooth family of symplectic forms $\om_{b,t}$
of the form $g_{b,t}^* \psi_b^*(\om_t)$, where $\psi_b$ is a symplectomorphism $(M_b,
\om_b) \to (M,\om)$.  Hence, for each $t$, $\ga^*(P)$ can be symplectically
identified with $$
\bigcup_{s\in [0,1]} (\{s\}\times (M, g_{s,t}^*(\om_t))),
$$
where $g_{1,t}^*(\om_t) = \om_t$ and the $g_{s,t}$ lie in an arbitrarily small neighborhood
$U$ of the identity in $\Diff(M)$.   By Moser's homotopy argument, we can choose $U$ so
small that each $g_{1,t}$ is isotopic to the identity in the group $\Symp(M, \om_t)$.  This
proves $(i)$.\QED

\begin{cor}\label{cor:obv}  The pullback of a stable Hamiltonian
bundle is stable.\end{cor}
\proof{} 
Suppose that  $P\to B$ is the pullback of $P'\to B'$
via $B\to B'$  so that there is a diagram
$$
\begin{array}{ccc}
P& \to& P'\\
\downarrow &&\downarrow\\
B&\to&B'.
\end{array}
$$
By hypothesis, the restriction $ H^2(P')\to H^2(M)$ is surjective. 
But this map factors as $H^2(P')\to H^2(P)\to H^2(M)$.  Hence
$H^2(P)\to H^2(M)$ is also surjective. 
\QED

\begin{lemma}\label{le:st1}\begin{itemize}
\item[(i)]
Every Hamiltonian bundle over $S^2$ is stable.
\item[(ii)]  Every symplectic bundle over a $2$-connected
base $B$ is Hamiltonian stable.
\end{itemize}
\end{lemma}
\proof{} 
$(i)$ holds because every Hamiltonian bundle over $S^2$ is c-split, in particular
the restriction map $H^2(P)\to H^2(M)$ is surjective .
$(ii)$ follows immediately from Lemma~\ref{le:ham2c}.\QED

\NI
{\bf Proof of Theorem~\ref{thm:stable}.}

This states that every Hamiltonian bundle is stable.
To prove this, first observe that we can restrict to the case when $B$ is
simply connected.
 For the map $B \to B\,\Ham(M)$ classifying $P$ factors through a map $C\to
B\,\Ham(M)$, where $C = B/B_1$ as before, and the stability of the induced
bundle over $C$ implies that for the original bundle by
Corollary~\ref{cor:obv}.

Next observe that by Lemma~\ref{le:st2}  a
Hamiltonian bundle $P\to B$ is stable if and only if the differentials
$d_k^{0,2}: E_k^{0,2}\to E_k^{k,3-k}$ in its Leray cohomology spectral sequence vanish  
on the whole of $H^2(M)$  for $k = 2,3$.  Exactly as in the proof of
Lemma~\ref{le:shamchar}  we can reduce the statement for $d_2^{0,2}$ to the
case $B = S^2$.  Thus $d_2^{0,2}= 0$ by Lemma~\ref{le:st1}$(i)$.  Similarly, we
can reduce the statement for $d_3^{0,2}$ to the case $B
= S^3$ and then use Lemma~\ref{le:st1}$(ii)$.\QED

\subsection{Functorial properties}\label{ss:funct}

We begin with some trivial observations and then discuss composites of Hamiltonian
bundles. The first lemma is true for any class of bundles with specified
 structural group.

 \begin{lemma}\label{le:pullb} Suppose that $\pi:P\to B$ is Hamiltonian and that $g:B'\to
B$ is a continuous map.  Then the induced bundle $\pi':g^*(P)\to B'$
is Hamiltonian.
\end{lemma}

Recall from \S\ref{ss:char} that any extension $\tau$ of the forms on the
fibers is called a connection form.

\begin{lemma}\label{le:sym} If $P\to B$ is a smooth Hamiltonian fiber bundle over a
symplectic base $(B,\si)$ and if $P$ is compact then  there is a connection form 
$\Om^\ka$ on $P$ that is
symplectic. \end{lemma}
\proof{}  By Proposition~\ref{prop:hamchar2}, the bundle $P$ 
carries a closed connection form $\tau$.
Since $P$ is compact, the form $\Om^\ka = \tau+\ka\pi^*(\si)$ is symplectic for large
$\ka$. \QED

Observe that the deformation type of the form $\Om^\ka$ is unique for sufficiently
large $\ka$ since given any two closed  connection forms $\tau,\tau'$   the linear isotopy
$$
t \tau + (1-t)\tau' + \ka\pi^*(\si),\quad 0\le t\le 1,
$$
consists of symplectic forms for sufficiently
large $\ka$.  However, it can happen that there is a symplectic connection
form $\tau$ such that $\tau + \ka\pi^*(\si)$ is not symplectic for small $\ka > 0$, even
though it is symplectic for large $\ka$.  (For example, suppose $P = M\times B$
and that $\tau$ is the sum $\om + \pi^*(\om_B)$ where $\om_B + \si$ is not
symplectic.)

Let us now consider the behavior of Hamiltonian bundles  under
composition. If  
$$
(M,\om)\to P
\stackrel{\pi_P}\to
X,\quad\mbox{and}\quad (F,\si)\to X\stackrel{\pi_X}{\to} B
$$
 are Hamiltonian
fiber bundles,  then the restriction 
$$
\pi_P:\quad W = \pi_P^{-1}(F)\longrightarrow F
$$
is a Hamiltonian fiber bundle.  Since $F$ is a manifold, we can assume
without loss of generality that $W\to F$ is smooth: see
Lemma~\ref{le:smooth}.  Moreover, since $(F, \si)$ is symplectic 
  Lemma~\ref{le:sym} implies that the manifold $W$ carries a
symplectic connection form $\Om_W^\ka$, and
 it is natural to ask when the composite
map $\pi:P\to B$ with fiber $(W, \Om_W^\ka)$  is itself Hamiltonian.   


\begin{lemma}\label{le:comp}  Suppose in the above situation
 that $B$ is  simply connected  and that $P$ is compact.
 Then  $\pi = \pi_X\circ \pi_P: P\to B$ is a 
Hamiltonian fiber bundle with fiber $(W, \Om_W^\ka)$, where
$ \Om_W^\ka = \tau_W  + \ka\pi_P^*(\si)$,
$\tau_W$ is any symplectic connection form  on $W$, 
and $\ka$ is sufficiently large.
 \end{lemma} 
\proof{} 
 By Lemma~\ref{le:smooth},  we may assume that the base 
$B$ as well as the fibrations
are smooth.  We first show that there is some symplectic form on $W$ for which
$\pi$ is Hamiltonian and then show that it is Hamiltonian with respect to 
the given form $\Om_W^\ka$.

Let $\tau_P$ (resp.
$\tau_X$) be a closed connection form for
the bundle $\pi_P$, (resp. $ \pi_X$),
and let $\tau_W$ be its restriction to $W$.
Then $\Om_W^\ka$ is the restriction to $W$ of the closed form
$$
\Om_P^\ka = \tau_P + \ka\pi_P^*(\tau_X).
$$
By increasing $\ka$ if necessary we can ensure that $\Om_P^\ka$ restricts to a symplectic
form on every fiber of $\pi$ not just on the the chosen fiber $W$.
 This shows firstly that
$\pi:P\to B$ is symplectic, because there is a well defined 
 symplectic form on each of its fibers, and secondly that it is Hamiltonian with respect to
this form $\Om_W^\ka$ on the fiber $W$.   Hence 
Lemma~\ref{le:st2} implies that  $H^2(P)$ surjects onto $H^2(W)$. 

Now suppose that  $\tau_W$ is {\it any}
closed connection form on $\pi_P: W\to F$.  
Because the restriction map $H^2(P) \to H^2(W)$ is surjective,
the cohomology class $[\tau_W]$ is the restriction of a class on $P$
and so, by  Thurston's construction, the form $\tau_W$ can be extended to a closed connection
form
$\tau_P$  for the bundle $\pi_P$. 
Therefore the previous argument applies in this case too. \QED

Now let us consider the general situation, when $\pi_1(B) \ne 0$.  
The proof of the lemma above applies to show that the composite bundle $\pi:P\to B$ 
is symplectic with respect to suitable $\Om_W^\ka$ and that it has a symplectic connection
form.  However,
 even though $\pi_X:X\to B$ is symplectically trivial over the $1$-skeleton 
of $B$ the
same may not be true of the composite map $\pi:P\to B.$  Moreover, in general it is not
clear whether triviality with respect to one form $\Om_W^\ka$ implies 
it for another.
Therefore, we may conclude the following:

\begin{prop}\label{prop:comp}
If
$(M,\om)\to P
\stackrel{\pi_P}\to
X,$ and $ (F,\si)\to X\stackrel{\pi_X}{\to} B$ are Hamiltonian fiber bundles and $P$ is
compact, then the composite
 $\pi =\pi_X\circ\pi_P:P\to B$ is a symplectic fiber bundle with respect to any form
$\Om_W^\ka$ on its fiber $W = \pi^{-1}(pt)$, provided that $\ka$ is sufficiently large. 
Moreover if $\pi$ is symplectically trivial over the $1$-skeleton of $B$ with respect to
 $\Om_W^\ka$ then $\pi$ is Hamiltonian.
\end{prop}

In practice, we will apply these results in cases 
where $\pi_1(B) = 0$.  We will not specify
the precise form on $W$, assuming that 
it is $\Om_W^\ka$ for a suitable $\ka$.

\section{Splitting of rational cohomology} \label{se:c-splitting}

We write $H_*(X), H^*(X)$ for the rational (co)homology of $X$.
Recall that a bundle $\pi:P\to B$ with fiber $M$ 
is said to be {\it c-split} if 
$$
H^*(P) \cong H^*(B)\otimes H^*(M).
$$
This happens if and only if $H_*(M)$ injects into $H_*(P)$.
Dually, it happens if and only if the restriction map 
$H^*(P)\to H^*(M)$ is onto.  
Note also that a bundle $P\to B$ c-splits if and only
 if the $E_2$ term of its cohomology
spectral sequence is a product and all the 
differentials $d_k, k\ge 2,$ vanish.

In this section we prove all parts of Theorem~\ref{thm:split}.
 We begin by using topological arguments that are based on
the fact that bundles over $S^2$ are c-split.  This was proved
in~\cite{LMP,Mc} by geometric arguments using
Gromov--Witten invariants.  In \S\ref{ss:gw} we discuss the extent to which
these geometric arguments generalize. Finally in \S\ref{ss:fex} we 
 discuss c-splitting in a homotopy-theoretic
context.

\subsection{A topological discussion of c-splitting}

The first lemma is obvious but useful.  We will often refer to its second part as the
Surjection Lemma. 

\begin{lemma}\label{le:funct} Consider a commutative diagram
$$
\begin{array}{ccc}
P' & \to & P\\
\downarrow & &\downarrow\\
B'& \to & B
\end{array}
$$
where $P'$ is the induced bundle.  Then:
\begin{itemize}
\item[(i)]  If $P\to B$ is c-split so is $P'\to B'$.\SmS

\item[(ii)] {\bf (Surjection Lemma)} If $P'\to B'$ is c-split and $H_*(B')\to
H_*(B)$ is surjective, then $P\to B$ is c-split.\end{itemize}
\end{lemma}
\proof{} $(i)$:  Use the fact that $P\to B$ is c-split if and only if the map $H_*(M)\to
H_*(P)$ is injective.

\NI
$(ii)$:  The induced map on the $E_2$-term of the cohomology spectral sequences is
injective.  Therefore the existence of a nonzero differential in the spectral sequence $P\to
B$ implies one for the pullback bundle $P'\to B'$. \QED

\begin{cor}\label{cor:bl}  Suppose that $P\to W$ is a Hamiltonian fiber bundle over
a symplectic manifold $W$ and that its pullback to some blowup $\widehat W$ of $W$
is c-split.  Then $P\to W$ is c-split.
\end{cor}
\proof{}  This follows immediately from $(ii)$ above since the map
$H_*(\widehat W) \to H_*(W)$ is surjective.\QED

\begin{lemma}\label{le:funct2}
If $(M,\om)\stackrel{\pi}\to  P\to B$ is a 
compact Hamiltonian  bundle over a simply connected CW-complex
$B$ and if every Hamiltonian fiber bundle over 
$M$ and $B$ is c-split, then every Hamiltonian
bundle over $P$ is c-split. 
\end{lemma}
\proof{}  Let $\pi_E: E\to P$ be a Hamiltonian bundle with fiber $F$ and let
$$
F\to W\to M
$$
be its restriction over $M$.    Then by assumption the latter bundle
c-splits so that $H_*(F)$ injects into $H_*(W)$.  Lemma~\ref{le:comp} 
implies that
the composite bundle $E\to B$ is
Hamiltonian with fiber $W$ and therefore also c-splits.   Hence $H_*(W)$
injects into $H_*(E)$. Thus $H_*(F)$ injects into $H_*(E)$, as required.\QED

\begin{lemma} \label{le:sur} If $\Si$ is a
closed orientable surface then any Hamiltonian bundle over 
$S^2 \times \ldots \times S^2 \times \Si$
 is c-split. 
\end{lemma}  
\proof{} Consider any degree one map $f$ from $\Si\to S^2$. 
Because $\Ham(M, \om)$ is
connected, $B\,\Ham(M, \om)$ is simply connected, and therefore any
homotopy class of maps from $\Si\to B\,\Ham(M, \om)$ factors through
$f$.  Thus any Hamiltonian bundle over $\Si$ is the pullback by $f$ of a Hamiltonian
bundle over $S^2$. Because such bundles  c-split
 over $S^2$, the same is true over $\Si$ by Lemma~\ref{le:funct}$(i)$.

The statement for $S^2 \times \ldots \times S^2 \times \Si$
 is now a direct consequence of iterative 
applications of Lemma~\ref{le:funct2} applied to the trivial bundles 
$
 S^2 \times \ldots \times S^2 \times \Si \to S^2.
$
\QED

\begin{cor} \label{cor:odd} Any Hamiltonian bundle over
$S^2 \times \ldots \times S^2 \times S^1$ is c-split.
\end{cor}

\begin{prop}\label{prop:sph} For each $k\ge 1$, every 
Hamiltonian bundle over $S^k$  c-splits.
\end{prop}
\proof{}  By Lemma~\ref{le:sur} and 
Corollary~\ref{cor:odd} there is for each
$k$ a  $k$-dimensional closed manifold 
$X$ such that every Hamiltonian bundle 
over $X$ c-splits.  Given any Hamiltonian
bundle $P\to S^k$ consider its pullback to 
$X$ by a map $f:X\to S^k$ of degree $1$.
Since the pullback c-splits, the original bundle does too by
the surjection lemma.  \QED

 As we shall see this result implies that the action of 
the homology groups of $\Ham(M)$ on $H_*(M)$ is always trivial.  
Here are some other examples of situations in which Hamiltonian
bundles are c-split.

\begin{lemma} \label{le:products} Every Hamiltonian 
bundle over $\CP^{n_1} \times \ldots \times \CP^{n_k}$ c-splits.
\end{lemma}
\proof{}
Let us prove first that it splits over $\CP^n$. 
Use induction over $n$.  Again it
holds when $n=1$. Assuming the result for
$n$ let us prove it for $n+1$. Let $B$  be the blowup of 
$\CP^{n+1}$ at one point. Then $B$
fibers over $\CP^n$ with fiber $\CP^1$. Thus every Hamiltonian
bundle over $B$ c-splits by Lemma~\ref{le:funct2}. The result for $\CP^{n+1}$ now follows
from Corollary~\ref{cor:bl}. Finally Hamiltonian bundles  c-split over products of 
projective spaces by repeated applications of Lemma~\ref{le:funct2}. \QED

\begin{cor}\label{cor:bt}  Every Hamiltonian bundle whose structural group
reduces to a subtorus $T \subset \Ham(M)$ c-splits.
\end{cor}
\proof{} It suffices to consider the universal model
$$
M \longrightarrow ET\times_T M \longrightarrow BT,
$$
and hence to show that all Hamiltonian bundles over $BT$ are c-split.
But this is equal to $\CP^{\infty} \times
\ldots \times \CP^{\infty}$ and the proof that the $i^{\rm th}$ group
of homology of the fiber injects in $P \to \CP^{\infty} \times
\ldots \times \CP^{\infty}$ may be reduced to the proof that it injects
in the restriction of the bundle $P$ over 
$\CP^j \times \ldots \times \CP^j$ for a sufficiently large
$j$. But this is Lemma~\ref{le:products}.\QED

\begin{rmk}\label{rmk:bt}\rm  Observe that the proof of the above corollary
shows that every Hamiltonian bundle over $\CP^{\infty} \times
\ldots \times \CP^{\infty}$ c-splits.  Since the structural group of
 such a bundle can be larger than the torus $T$, our result extends
the  Atiyah-Bott splitting theorem for Hamiltonian
bundles with structural group $T$.
\end{rmk}

For completeness, we show how the above corollary leads to a proof of the 
splitting of $G$-equivariant cohomology where $G$ is a Lie subgroup of 
$\Ham(M, \om)$.  

\begin{cor}\label{cor:funct2}
If $G$ is a compact connected Lie group that acts in a Hamiltonian way on $M$
then any bundle $P\to B$ with fiber $M$ and structural group $G$ is c-split. 
In particular, 
$$
H^*_G(M) \cong H^*(M)\otimes H^*(BG).
$$  
\end{cor}

\proof{}  By Lemma~\ref{le:funct}$(i)$ we only need to prove the second
statement, since 
$$
M_G = EG\times _GM \longrightarrow  BG
$$
is the universal bundle.
Every compact connected Lie group $G$ is the image of a homomorphism
$T\times H \to G$, where the torus $T$ maps onto the identity component of the
center of $G$ and $H$ is the semi-simple Lie group corresponding to the
commutator subalgebra $[{\rm Lie} (G), {\rm Lie} (G)]$ in the Lie algebra
${\rm Lie}(G)$.  It is easy to see that this  homomorphism induces a
surjection on rational homology $BT \times BH \to BG$.  Therefore, we may
suppose that $G = T\times H$.  
Let $T_{max} = (S^1)^k$ be the
maximal torus of the semi-simple group $H$. Then 
the induced map on cohomology $H^*(BH) \to H^*(BT_{max}) =
\Q[a_1,\dots, a_k]$ takes $H^*(BH)$ bijectively onto the set of polynomials in
$H^*(BT_{max})$ that are invariant under the action of the Weyl group, by the
Borel-Hirzebruch theorem. Hence  the maps $BT_{max} \to BH$
and $BT \times BT_{max} \to BG$
induce a surjection on homology.   Therefore the desired result
follows from the surjection lemma and the last corollary (or the 
Atiyah--Bott theorem.) 
  \QED

\begin{lemma}\label{le:bott}  Every
Hamiltonian bundle over a coadjoint orbit c-splits.
\end{lemma}
\proof{}  This is an immediate
consequence of the results by Grossberg--Karshon~\cite{GK}\S3 on Bott towers.  A
Bott tower is an iterated fibration $M_k\to M_{k-1} \to\dots \to M_1 = S^2$ of
K\"ahler manifolds where each map $M_{i+1}\to M_i$ is a fibration with fiber
$S^2$.  They show that any coadjoint orbit $X$ can be blown up to a
manifold that is diffeomorphic to a Bott tower $M_k$.  Moreover the
blow-down map $M_k\to X$ induces a surjection on
rational homology.   Every Hamiltonian bundle over $M_k$ c-splits by repeated
applications of Lemma~\ref{le:funct2}. Hence the result follows from the
surjection lemma.\QED

 \begin{lemma}\label{le:3c} Every Hamiltonian bundle over a
$3$-complex $X$ c-splits. \end{lemma}
\proof{}  As in the proof of stability given in Theorem~\ref{thm:stable} 
we can reduce this to the cases $X= S^2$ and $X= S^3$ and then use
Proposition~\ref{prop:sph}.  The only difference from 
the stability result is that we now
require the differentials $d_2^{0,q}, d_3^{0,q}$ to 
vanish for all $q$ rather than just at
$q=2$.\QED

 \begin{lemma}\label{le:spheres} Every Hamiltonian bundle over
a product of spheres c-splits, provided that
there are no more than $3$ copies of $S^1$. \end{lemma}
\proof{}  By hypothesis $B = \prod_{i\in I} S^{2m_i}\times \prod_{j\in
J}S^{2n_i + 1}\times T^k$,  where $n_i > 0$ and $0\le k\le 3$.  Set
 $$
B' =  \prod_{i\in I} \CP^{m_i} \times \prod_{j\in J} {\CP}^{n_i} \times T^{|J|}
\times T^\ell,
$$
where $\ell = k$ if $k + |J|$ is even and $= k+1$ otherwise.
Since $\CP^{n_i}\times S^1$ maps onto $S^{2n_i+1}$ by a map of degree $1$,
there is a homology surjection $B'\to B$ that maps the factor $T^\ell$ to
$T^k$. By the surjection lemma, it suffices to show  that  the pullback bundle
$P'\to B'$ is c-split. 

Consider the fibration 
$$
T^{|J|}\times T^\ell \to B' \to 
\prod_{i\in I} \CP^{m_i} \times \prod_{j\in J} {\CP}^{n_i}.
$$
Since $|J| + \ell$ is even, we can think of this as a Hamiltonian bundle.
Moreover, by construction, the restriction of the bundle $P'\to B'$ to
$T^{|J|}\times T^\ell$ is the pullback of a bundle over $T^k$, since the map
$T^{|J|}\to B$ is nullhomotopic. (Note that each $S^1$ factor in  $T^{|J|}$ goes
into a different sphere in $B$.)   Because $k\le 3$, the bundle over $T^k$ c-splits.
Hence we can apply the argument in Lemma~\ref{le:funct2} to conclude that 
$P'\to B'$ c-splits.\QED

\begin{lemma}\label{le:dim2}
 Every Hamiltonian bundle  whose fiber has 
cohomology  generated by $H^2$ is c-split.
\end{lemma}
\proof{} This is an immediate consequence of 
Theorem~\ref{thm:stable}.\QED

\NI
{\bf Proof of Theorem~\ref{thm:split} }

Parts $(i)$ and $(iii)$ are proved in the
Lemmas~\ref{le:bott},~\ref{le:3c} and~\ref{le:spheres} above,
and part $(ii)$ is proved in~\S\ref{ss:gw} below.\QED

Of course, the results in this section can be extended further
by applying the surjection lemma and
variants of Lemma~\ref{le:funct2}.   For example, any
Hamiltonian fibration c-splits if its  base  $B$ is the image of a homology
surjection from a product of spheres and projective spaces, provided that
there are no more than three $S^1$ factors.  One can also consider iterated
fibrations of projective spaces, rather than simply products.  
However, we have not yet managed to deal with arbitrary products of
spheres.   In order to do this, it would
suffice to show that every Hamiltonian bundle over a torus $T^{m}$ c-splits. 
This question has not yet been resolved for $m\ge 4$.

\subsection{Hamiltonian bundles and Gromov--Witten invariants}
\label{ss:gw}

We begin by sketching an alternative proof that every Hamiltonian bundle over
$B = \CP^n$ is c-split that  generalizes the arguments in~\cite{Mc}.  We will
use the language of~\cite{Mcv}, which is  based on the Liu--Tian~\cite{LiuT}
approach to general Gromov--Witten invariants.
No doubt any treatment of
general Gromov--Witten invariants could  be used instead.\MS

\NI
{\bf Proof that every Hamiltonian bundle over $\CP^n$ is c-split.}

The basic idea is to show that 
the inclusion $\io: H_*(M) \to H_*(P)$ is injective by showing that 
for every nonzero $a\in H_*(M)$ there is $b\in H_*(M)$ and $\si\in H_2(P;\Z)$
for which the  Gromov--Witten invariant $n_P(\io(a), \io(b);\si)$
is nonzero.  Intuitively speaking this invariant
``counts the number of isolated $J$-holomorphic curves in $P$ that represent
the class $\si$ and meet the  classes $\io(a), \io(b)$."  More correctly, it is
defined to be the intersection number of the image of the evaluation map
$$
ev:\oMm\,\!^\nu_{0,2}(P,J,\si)\longrightarrow  P\times P
$$
with the class $\io(a)\times \io(b)$, where
$\oMm\,\!^\nu_{0,2}(P,J,\si)$  is a virtual moduli cycle
made from perturbed $J$-holomorphic curves with $2$ marked points,
and $ev$ is given by evaluating at these two points.
As explained in~\cite{Mcv,Mc}, $\oMmn = \oMm\,\!^\nu_{0,2}(P,J,\si)$
is a branched pseudomanifold, i.e. a kind of stratified space whose top
dimensional strata are oriented and have rational labels.  Roughly speaking, one
can think of it as a finite simplicial complex, whose dimension $d$ equals the 
``formal dimension" of the moduli space, i.e. the index of the relevant operator.
The  elements of $\oMmn$ are stable maps $[\Si, h, z_1,z_2]$
where $z_1,z_2$ are two marked points on the nodal, genus $0$, Riemann
surface $\Si$, and the map $h:\Si\to P$  satisfies a perturbed Cauchy--Riemann
equation $\pJ h = \nu_h$.  The perturbation $\nu$ can be arbitrarily small, and
is chosen so that each stable map in $\oMmn$  is a regular point for the
appropriate Fredholm operator.  Hence $\oMmn $ is
often called a {\it regularization} of the unperturbed moduli space $\oMm=
\oMm_{0,2}(P,J,\si)$ of all $J$-holomorphic stable maps.

Given any Hamiltonian bundle $P_S\to S^2$ 
and any $a\in H_*(M)$, it was shown in \cite{LMP,Mc} that
there is $b\in H_*(M)$ and a lift
 $\si_S\in H_2(P_S;\Z)$ of the fundamental class of $S^2$ to $P_S$ such that
$$
n_{P_S}(\io_S(a), \io_S(b); \si_S) \ne 0,
$$
where $\io_S$ denotes the inclusion into $P_S$.
Therefore, if $P_S$ is identified with the restriction of $P$ to a complex line
$L_0$  in $B$ and if $a,b$ and $\si_S$ are as above, it
suffices to prove that 
$$
 n_{P_S}(\io_S(a), \io_S(b); \si_S) = n_P(\io(a), \io(b);\si)
$$
where $\si$ is the image of $\si_S$ in $P$. Note that a direct count shows that
the dimensions of the appropriate virtual moluli spaces
$\oMm\,\!^\nu_{0,2}(P_S,J_S,\si_S)$ and
$\oMm\,\!^\nu_{0,2}(P,J,\si)$ differ by the codimension of
$P_S \times P_S$ in $P \times P$ (which equals the codimension of
$\CP^1 \times \CP^1$ in $\CP^n \times \CP^n$) so that the both sides
are well-defined.

As was shown in~\cite{Mc}~Corollary~4.11,  one can construct the virtual
moduli cycle
$\oMmn(P_S) = \oMm\,\!^\nu_{0,2}(P_S,J_S,\si_S)$ using an almost complex
structure $J_S$ and a perturbation $\nu$ that are compatible with the bundle. 
In particular, this implies  that the projection $P_S\to S^2$ is $(J_S,
j)$-holomorphic (where $j$ is the usual complex structure on $S^2$) and that
every element of $\oMmn(P_S)$ projects to a
$j$-holomorphic stable map in $S^2$.

We claim that this is also true for the bundle $P\to B$.
In other words, we can choose $J$ so that the projection $(P,J) \to (B,j)$ is
holomorphic, where $j$ is the usual complex structure on $B = \CP^n$, and we can
choose $\nu$ so that every element in $\oMmn(P)$
projects to a $j$-holomorphic  stable map in $B$.
The proof is
exactly as before: see~\cite{Mc} Lemma 4.9.  The essential point  is
that every element of the unperturbed moduli space 
$\oMm_{0,2}(\CP^n, j, [\CP^1])$ is regular.
In fact, the top stratum in
$\oMm_{0,2}(\CP^n, j, [\CP^1])$ is the space $\Ll = \Mm_{0,2}(\CP^n, j, [\CP^1])$ 
of all   lines in $\CP^n$ with $2$ distinct marked points.  The other stratum
completes this space by adding in the
lines with two coincident marked points, which are represented as
stable maps by a line together with a ghost bubble
containing the two points.  

It follows that there is a projection map 
$$
proj:\;\; \oMm\,\!^\nu_{0,2}(P,J,\si) \longrightarrow \oMm_{0,2}(\CP^n, j,
[\CP^1]).
$$
 Moreover the inverse image of a line $L\in \Ll$
can morally speaking be identified with
$\oMm\,\!^\nu_{0,2}(P_S,J_S,\si_S)$.   The latter statement would be  correct
 if we were considering ordinary moduli spaces of stable maps, but the
virtual moduli space is not usually built in such a way that the fibers $
(proj)^{-1}(L)$ have the needed
structure of a branched pseudomanifold.  However, we can choose to
construct $\oMm\,\!^\nu_{0,2}(P,J,\si)$ 
so that this is true for all lines near a fixed
line $L_0$.  
In~\cite{Mc} (see also~\cite{Mcv}) a detailed recipe is given for
constructing  $\oMmn$ from the unperturbed moduli space $\oMm$. The
construction is based on the choice of suitable covers $\{U_i\}, \{ V_I\}$ of
$\oMm$ and of perturbations $\nu_i$ over each $U_i$.
Because regularity is an open condition, one can make these choices first for all
stable maps that project to the fixed line $L_0$ and then extend
to the set of stable maps that project to nearby lines in such a way that 
$\oMmn$  is locally a product near the  fiber over $L_0$:
see the proof of \cite{Mc}~Proposition~4.6 for a very similar construction.

Once this is done, the rest of the argument  is easy.  If we identify $P_S$ with
$\pi^{-1}(L_0)$ and choose a representative $\al\times \be$ of $\io_S(a) \times
\io_S(b)$ in $P_S\times P_S$ that is transverse to the evaluation map from 
$\oMm\,\!^\nu_{0,2}(P_S,J_S,\si_S)$,  its image in $P\times P$ will be transverse
to the evaluation map from $\oMm\,\!^\nu_{0,2}(P,J,\si)$ because $proj$ is a
submersion at $L_0$.   Moreover, by \cite{Mc}~Lemma~4.14, we may suppose
that  $\al$ and $\be$ lie in distinct fibers of the projection $P_S\to S^2$.  Let
$x_a, x_b$ be the corresponding points of $\CP^n$ under the identification
$S^2 = L_0$.  Then  every stable map that contributes to $ n_P(\io(a), \io(b);\si)$
projects to an element of $\oMm_{0,2}(\CP^n, j, [\CP^1])$ whose marked points
map to the distinct points $x_a, x_b$.  Since there is a unique line in 
$\CP^n$ through two given points, in this case $L_0$, every stable map
that contributes to $ n_P(\io(a), \io(b);\si)$ must project to $L_0$ and hence 
be contained in 
$\oMm\,\!^\nu_{0,2}(P_S,J_S,\si_S)$. One can then  check that
$$
 n_{P_S}(\io_S(a), \io_S(b); \si_S) = n_P(\io(a), \io(b);\si),
$$
as claimed. The only delicate point here is the verify that the sign of each
stable map on the left hand side is the same as the sign of the corresponding
map on the right hand side. But this is also a consequence of the local
triviality of the above projection map (see \cite{L} for more details).\QED

The above argument generalizes easily to the case when $B$ is a
complex blowup of $\CP^n$.

\begin{prop}\label{prop:blcp}  Let $B$ be a blowup of $\CP^n$ 
along a disjoint
union $Q = \coprod Q_i$ of complex submanifolds, 
each of complex codimension $\ge 2$. 
Then every Hamiltonian bundle $B$  is c-split. \end{prop}
\proof{}  
 The above argument applies almost verbatim in the case when $Q$ 
is a finite set of points.  The top stratum of 
 $\oMm_{0,2}(B, j, [\CP^1])$ still consists of lines marked by two distinct points,
and again all elements of this unperturbed moduli space are regular.

 In the general case, there is a blow-down map $\psi: B\to \CP^n$ which is
bijective over $\CP^n - Q$, and we can choose $j$ on $B$ so that
 the exceptional divisors $\psi^{-1}(Q)$ are $j$-holomorphic, and so
that $j$ is pulled back from the usual structure on $\CP^n$ outside a small 
neighborhood of  $\psi^{-1}(Q)$.  Let $L_0$ be a complex line in $\CP^n - Q$.
Then its pullback to $B$ is still $j$-holomorphic.   Hence, although the
unperturbed moduli space  $\oMm_{0,2}(B, j, [\CP^1])$ may contain nonregular
and hence  ``bad" elements, its top stratum does contain an open set
$U_{L_0}$ consisting of marked lines near $L_0$ that are regular. 
Moreover, if we fix two points $x_a, x_b$ on $L_0$, every element of
$\oMm_{0,2}(B, j, [\CP^1])$ whose marked points map sufficiently close to
$x_a, x_b$ actually lies in this open set $U_{L_0}$.
 We can
then regularize  $\oMm_{0,2}(B, j, [\CP^1])$ to a virtual moduli cycle that
contains the open set $U_{L_0}$ as part of its top stratum.  Moreover, because
the construction of the regularization is local with respect to 
$\oMm_{0,2}(B, j, [\CP^1])$,  this
regularization  $\oMm\,\!^\nu_{0,2}(B, j, [\CP^1])$ 
will still have the property that 
each of its elements whose marked points map sufficiently close to
$x_a, x_b$ actually lies in this open set $U_{L_0}$.

We can now carry out the previous argument, choosing $J$ on $P$ to
be fibered, and constructing $\nu$ to be compatible with the fibration
on that part of $\oMm_{0,2}(P, J, \si)$ that projects to $U_{L_0}$.
Further details will be left to the reader.\QED

\begin{cor}\label{cor:bl1}  Let $X = \# k\CP^2\# \ell \bcp$ be the connected
sum of  $k$ copies of $\CP^2$ with $\ell$ copies of $\bcp$.  If one of $k, \ell$
is $\le 1$ then every bundle over $X$ is c-split.
\end{cor}
\proof{}  By reversing the orientation of $X$ we can suppose that $k\le 1$.
The case $k=1$ is covered in the previous proposition.  When $k=0$,
 pull the bundle back over the blowup of $X$ at one
point and then use the Surjection Lemma~\ref{le:funct}(ii).\QED

\MS

 The previous proof can easily be generalized to the case
of a symplectic base $B$ that has a spherical 2-class $A$
with Gromov-Witten invariant of the form $n_B(pt, pt, c_1,\dots, c_k; A)$ 
absolutely equal
to $1$. (By this we mean that for some  generic $j$ on $B$
the relevant moduli space contains exactly one element,
which moreover parametrizes an embedded curve in $B$.)

Here $c_1,\dots, c_k$
are arbitrary homology classes of $B$ and we assume that $k \ge 0$.
Again the idea is to construct  the regularizations
$\oMm\,\!^\nu_{0,2+k}(P,J,\si)$ 
and $\oMm\,\!^\nu_{0,2+k}(B,j,A)$ so that there is
a  projection from one to the other which is a fibration at least near the
element of $\oMm\,\!^\nu_{0,2+k}(B,j,A)$ that contributes to $n_B(pt, pt,
c_1,\dots, c_k; A)$.   Thus $B$ could be the blowup of $\CP^n$ along a
symplectic submanifold $Q$ that is disjoint from a complex line.
One could also take similar blowups of products of projective spaces, or, more
generally, of iterated fibrations of projective spaces.  For example, if 
$B = \CP^m\times
\CP^n$ with the standard complex structure then there is one complex line in
the diagonal class $[\CP^1] + [\CP^1]$ passing through any two points and
a cycle  $H_1 \times H_2$, where $H_i$ is the hyperplane
class, and one could blow up
along any symplectic submanifold that did not meet one such line.

 It is also very likely that this argument 
can be extended to apply when
all we know about  $B$ is that some
Gromov--Witten invariant  $n_B(pt, pt, c_1,\dots, c_k; A)$ 
is nonzero, for example, if $B$ is a blowup of $\CP^n$ along 
arbitrary symplectic
submanifolds.  There are two new problems here: $(a)$  we must control the
construction of $\oMm\,\!^\nu_{0,2+k}(P,J,\si)$ in a neighborhood of all the
curves that contribute to $n_B(pt, pt, c_1,\dots, c_k; A)$ and $(b)$  we must
make sure that the orientations are compatible so that
 curves in $P$  projecting over different and noncancelling curves in $B$ do
not cancel each other in the global count of the Gromov-Witten invariant in
$P$.   Note that the bundles given by restricting $P$ to the different
curves counted in  $n_B(pt, pt, c_1,\dots, c_k; A)$ are diffeomorphic, 
since, this being a homotopy theoretic question, we can always replace $X$ by the
simply connected space $X/(X_1)$ in which these curves are homotopic: see
Corollary~\ref{cor:hamchar2}. Thus what is needed for this generalization
to hold is to develop further the theory of fibered GW-invariants that
was begun in~\cite{Mc}. See \cite{L}.

\subsection{Homotopy-theoretic reasons for c-splitting} \label{ss:fex}

In this section we discuss c-splitting  in a
homotopy-theoretic context.  Recall
that a c-Hamiltonian bundle is a    smooth bundle
$P\to B$ together with a class $a\in H^2(P)$ whose restriction $a_M$ to the
fiber  $M$ is c-symplectic, i.e. $(a_M)^n\ne 0$ where $2n = \dim(M)$.
Further  a closed manifold $M$ is said to
satisfy the hard Lefschetz condition with respect to the class
$a_M\in H^2(M, \R)$ if the maps
$$
\cup (a_M)^k: H^{n-k}(M,\R) \to H^{n+k}(M,\R),\quad 1\le k\le n,
$$ 
are isomorphisms.   In this case, elements in $H^{n-k}(M)$ that vanish when
cupped with $(a_M)^{k+1}$ are called primitive, and the cohomology of
$M$ has an additive basis consisting of elements of the form $b\cup (a_M)^\ell$
where $b$ is primitive and $\ell \ge 0$.  (These 
manifolds are sometimes called \lq\lq cohomologically K\"ahler.")

\begin{lemma}\label{le:blan}{\bf (Blanchard~\cite{Bl})}   Let  $M\to P\to B$
be a c-Hamiltonian bundle  such that $\pi_1(B)$ acts trivially on $H^*(M,
\R)$.  If in addition $M$ satisfies the hard Lefschetz condition with respect to
the c-symplectic class $a_M$, then the bundle c-splits. \end{lemma}
\proof{}
The proof is by contradiction.  Consider the Leray spectral sequence in 
cohomology  and suppose that
$d_p$ is the first non zero differential.  Then, $p\ge 2$ and
 the $E_p$ term in the spectral sequence is isomorphic to
the $E_2$ term and so can be identified with the tensor product $H^*(B)\otimes
H^*(M)$. Because of
the product structure on the spectral sequence, one of the differentials $d_p^{0,i}
$ must be nonzero.  So there is $b\in E_p^{0,i} \cong H^i(M)$ such that
$d_p^{0,i}(b) \ne 0$.  We may assume that $b$ is primitive (since these
elements together with $a_M$ generate $H^*(M)$.)
Then $b \cup a_M^{n-i} \ne 0$ but $b \cup a_M^{n-i+1} = 0$.
 
We can write $d_p(b) = \sum_j e_j\otimes f_j$ where $e_j \in H^*(B)$
and $f_j \in H^\ell(M)$ where $\ell < i$.  Hence $f_j \cup a_M^{n-i+1}\ne 0$
for all $j$ by the Lefschetz property.  Moreover, because the $E_p$ term is a
tensor product
 $$
(d_p(b))\cup a_M^{n-i+1} = \sum_j e_j \otimes (f_j \cup a_M^{n-i+1}) \ne 0.
$$
But this  is impossible since this element is the image via $d_p$ of the 
trivial element $b\cup a_M^{n-i+1}$.
\QED
 
Here is a related argument due to Kedra.\footnote
{Private communication}

\begin{lemma}\label{le:4f}
Every Hamiltonian bundle with $4$-dimensional fiber  c-splits.
\end{lemma}
\proof{}  Consider the spectral sequence as above. We know 
as in the proof of Lemma~\ref{le:3c} that $d_2 = 0$ and $d_3 = 0$.
  Consider $d_4$.  We just have
to check that $d_4^{0,3} = 0$ since $d_4^{0,i} = 0$ for $i = 1,2$ for dimensional
reasons, and $=0$ for $i = 4$ since the top class survives.
 
Suppose $d_4(b) \ne 0$ for some $b\in H^3(M)$.  Let $c\in H^1(M)$ be 
such that $b\cup c \ne 0$. Then $d_4(c) = 0$ and $d_4(b\cup c) = d_4(b)\cup c \ne
0$ since $d_4(b) \in H^4(B)\otimes H^0(M)$.  But we need $d_4(b\cup c) = 0$
since the top class survives.
So $d_4 = 0$ and then $d_k = 0$, $k> 4$ for reasons of dimension. \QED

Here is an example of a c-Hamiltonian bundle over
$S^2$ that is not c-split.  This shows that c-splitting is a geometric rather than
a topological (or homotopy-theoretic) property.
 
\MS \NI
{\bf Proof of  Lemma~\ref{le:cham}}

The idea is very simple.  First observe that if $S^1$ acts
on manifolds $X,Y$ with fixed points $p_X, p_Y$ then we can extend the $S^1$
action to the  connected sum $X\# Y^{opp}$  at $p_X, p_Y$ whenever the
 $S^1$ actions on the tangent spaces at $p_X$ and $p_Y$ are the same.  (Here
$Y^{opp}$ denotes $Y$ with the opposite orientation.) 
Now let $S^1$ act on $X = S^2\times S^2\times S^2$ by the diagonal action in
the first two spheres (and trivially on the third) and let the $S^1$ action on $Y$
be the example constructed in~\cite{Mcf} of a nonHamiltonian $S^1$ action
that has fixed points.  The fixed points in $Y$ form a disjoint union of $2$-tori
and the $S^1$ action in the normal directions has index $\pm 1$.  In other
words, there is a fixed point $p_Y$ in $Y$ at which we can identify $T_{p_Y}Y$
with $\C\oplus \C \oplus \C$, where $\theta\in S^1$ acts by multiplication by
$e^{i\theta}$ in the first factor, by multiplication by  $e^{-i\theta}$ in the second
and trivially in the third.  Since there is a fixed point on $X$ with the same
local structure, the connected sum $Z = X\# Y^{opp}$ does support an $S^1$-
action. Moreover $Z$ is a c-symplectic manifold.  There are many
possible choices of c-symplectic class:  under the
obvious identification of $H^2(Z)$ with $H^2(X) + H^2(Y)$ we will take the 
c-symplectic class on $Z$ to be given by the class of the symplectic form on
$X$.

Let $P_X\to S^2$, $P_Y\to S^2$ and $P_Z\to S^2$ be the corresponding
bundles.   Then $P_Z$  can be thought of as the connect sum of $P_X$ with
$P_Y$ along the sections corresponding to the fixed points.
By analyzing the corresponding Mayer--Vietoris sequence, 
it is easy to check that the c-symplectic class
on $Z$ extends to $P_Z$.  Further, the fact that the symplectic class in
$Y$ does not extend to $P_Y$ implies that it does not extend to $P_Z$
either.  Hence $P_Z\to S^2$ is not c-split.\QED

\section{Action of the homology of $\Ham(M)$ on $H_*(M)$}\label{ss:act}

The action  $\Ham(M)\times M \to M$ 
gives rise to maps  
$$
H_k(\Ham(M))\times H_*(M) \to H_{*+k}(M):\quad (\phi, Z)\mapsto \tr_\phi(Z).
$$
Theorem~\ref{thm:act} states:

\begin{prop}\label{prop:trivact}  These maps are  trivial when $k\ge 1$.
\end{prop}

\proof{} To see this, let us first consider the action of a spherical element
$$
\phi:S^n\to \Ham(M).
$$
It is not hard to check 
that the homomorphisms 
$$
\tr_\phi:  H_k(M) \to H_{k+n}(M)
$$
are precisely the connecting
homomorphisms in the Wang sequence of the bundle 
$P_\phi\to S^{k+1}$ with clutching
function $\phi$:  i.e. there is an exact sequence
$$
\dots H_k(M)\stackrel{\tr_\phi}\to H_{k+n}(M) \to  H_{k+n}(P) \stackrel{\cap [M]}\to
H_{k-1}(M) \to\dots
$$
Thus the fact that $P_\phi\to S^{k+1}$ is c-split immediately implies that the $\tr_\phi$
are trivial.

Next recall that in a $H$\/-space the rational cohomology ring is generated by 
elements dual to the rational
homotopy. It follows that there is a basis for $H_*(\Ham(M))$ that is represented by
cycles of the form 
$$
\phi_1\times\dots\times \phi_k: S_1 \times \ldots \times S_k \to \Ham(M),
$$
where the $S_j$s are spheres and one defines the product of maps by using the product
structure in $\Ham(M)$.  Therefore it suffices to show that these product  elements act
trivially.  However, 
if $a \in H_{*}(M)$ is represented by the cycle $\al$,  then $\tr_{S_k}(\al)$ is
null-homologous, and so equals the boundary $\partial \be $ of some chain
$\be$. Therefore: $$
\partial\left( \tr_{S_1 \times \ldots \times S_{k-1}}(\be)\right) = \pm \tr_{S_1
\times \ldots \times S_{k-1}}(\partial \be) = \pm \tr_{S_1 \times \ldots \times 
S_{k-1}}(\tr_{S_k}(\al))  = \tr_{S_1 \times \ldots  \times S_k}(\al).
$$
Hence $ \tr_{S_1 \times \ldots \times S_k}(a) = 0$.
This completes the proof. \QED

\begin{prop}\label{prop:equiv} Let $P  \to B$ be a trivial 
symplectic bundle.   Then any Hamiltonian automorphism
$\Phi\in \Ham(P,\pi)$ 
acts as the identity map on $H_*(P)$. 
\end{prop}
\proof{} An element $\Phi \in \Ham(P, \pi)$ is a map of the form
$$
\Phi : B \times M \to B \times M:\quad (b,x) \mapsto (b, \Phi_b(x))
$$
where $\Phi_b\in \Ham(M)$ for all $b\in B$.  Let us denote the induced map
$B\times M\to M: (b,x)\mapsto \Phi_b(x)$ by $\al_\Phi$.  
The previous proposition implies that if $B$ is a closed manifold of dimension $> 0$,
or, more generally, if it carries  a fundamental cycle $[B]$ of degree $> 0$,
 $$
(\al_\Phi)_*([B]\otimes m) = \tr_{[B]} (m) = 0,\quad \mbox{for all } m\in H_*(M).
$$
We can also think of $\Phi: B \times M \to B \times M$  as the
composite  
$$
B\times M \stackrel{diag_B\times id_M}\longrightarrow B \times B\times M 
\stackrel{id_B\times \al_\Phi}\longrightarrow B\times M.
$$
The diagonal class in $B\times B$ can be written  as
$[B]\otimes [pt] + \sum_{i\in I} b_i\otimes b_i'$ where 
$b_i, b'_i \in H_*(B)$ with $\dim(b'_i) > 0$.\footnote
{ This holds because the projection onto the first factor
takes the  diagonal class onto the fundamental class of $B$.  When the base
is a closed manifold, the diagonal is represented by $\sum_{i \in I}
(-1)^{\dim b_i}  b_i \otimes {b_i'}$ where $\{b_i\}$ is a basis for $H_*(B)$
and $\{b_i'\}$ is its Poincar\'e dual.}  
Hence 
$$
\Phi_{\ast}([B] \otimes m) = [B]\otimes m +  \sum_{i \in I} 
                               b_i \otimes \tr_{b_i'}(m) = [B]\otimes m,
$$
where the last equality comes from Proposition~\ref{prop:trivact}.
More generally, given any class $b\in H_*(B)$,  represent it by the image of
the fundamental class $[X]$ of some polyhedron under a suitable map $X\to B$
and consider the pullback bundle $P_X\to X$.  
 Since the class $\Phi_*([X]\otimes
m)$   is represented by a cycle in $X\times M$ for any $m \in H_*(M)$, we can work
out what it is by looking at the pullback of $\Phi$ to $X\times M$. The argument
above then applies to show that $ \Phi_*([X]\otimes m) = [X]\otimes m$ whenever
$b$ has degree $> 0$. Thus $\Phi_* = id$ on all cycles in $H_{*>0}(B)\otimes
H_*(M)$.  However, it clearly acts as the identity on $H_0(B)\otimes H_*(M)$
since the restriction of $\Phi$ to any fiber is isotopic to the identity.\QED

A natural conjecture is that the analog of
Proposition~\ref{prop:equiv} holds for all Hamiltonian bundles.  We now show that
there is a close relation between this question  and the
problem of c-splitting of bundles. Given an automorphism $\Phi$ of a symplectic
bundle $M \to P\to B$ we define 
$P_{\Phi} = (P \times [0,1])/\Phi$ to be the corresponding bundle over  $B
\times S^1$.  If the original bundle and the automorphism are Hamiltonian, 
so is
$P_\Phi\to B\times S^1$, though the associated  bundle $P_\Phi\to B\times S^1 \to S^1$ 
over $S^1$ will not
be, except in the trivial case when $\Phi$ is in the identity component 
of $\Ham(P,\pi)$.

\begin{prop} \label{prop:aut-versus-c-split} Assume that  
a given Hamiltonian bundle $M \to P
\to B$ c-splits. Then a Hamiltonian automorphism $\Phi \in \Ham(P, \pi)$ acts trivially
(i.e. as the identity) on $H_*(P)$ if and only if the corresponding Hamiltonian
 bundle $P_{\Phi}\to B \times S^1$
c-splits. \end{prop}
\proof{} Clearly, the fibration $P\to B$ c-splits if and only if
every basis of the $\Q$\/-vector
space $H^*(M)$ can be extended to a set of classes in $H^*(P)$ that form a
basis  for a complement to the kernel  of the restriction map.  We will call
such a set of classes a
{\it Leray--Hirsch basis.}  It corresponds to a choice of splitting 
isomorphism $H^*(P) \cong H^*(B)\otimes H^*(M).$
Now, the only obstruction to extending a Leray--Hirsch
basis from $P$ to $P'$ is the nontriviality of the action of $\Phi$ 
on $H^*(P)$. This shows the ``only if '' part. 
 
Conversely, suppose that $P_\Phi$ c-splits and let $e_j, j\in J,$ be a Leray--Hirsch
basis for  $H^*(P_\Phi)$.   Then $H^*(P_\Phi)$ has a basis of the form 
$e_j\cup \pi^*(b_i), e_j\cup \pi^*(b_i \times [dt])$ where $b_i$ runs through a basis
for $H^*(B)$ and $[dt] $ generates $H^1(S^1)$.  Identify $P$ with $P\times \{0\}$
in $P_\Phi$ and consider some cycle $Z\in H_*(P)$.  Since the cycles $\Phi_*(Z)$ and
$Z$ are homologous in $P_\Phi$, the classes $e_j\cup \pi^*(b_i)$ have equal values on
$\Phi_*(Z)$ and $Z$.  But the restriction of these classes to $P$ forms a basis for
$H^*(P)$.  It follows that $[\Phi_*(Z)] = [Z]$ in $H_*(P)$.\QED

\begin{prop} \label{prop:automorphisms} Let $P\to B$ be a Hamiltonian
bundle.  Then the group $\Ham(P, \pi)$ acts trivially on
$H_*(P)$ if the base 
\begin{itemize}
\item[(i)] has dimension $\le 2$, or
\item[(ii)] is a product of spheres and projective spaces with no more than
two $S^1$ factors, or
\item[(iii)] is simply connected and has the property that all Hamiltonian
bundles over $B$ are c-split.
\end{itemize}
\end{prop}
\proof{} In all cases, the hypotheses imply that $P\to B$ c-splits. 
Therefore the previous proposition applies
 and $(i)$, $(ii)$ follow immediately from Theorem~\ref{thm:split}.
To prove $(iii)$, suppose that $B$ is a simply
connected compact CW complex over which every Hamiltonian fiber bundle c-splits.
Let $M \hookrightarrow P \to B \times S^1$ be any Hamiltonian bundle -- in
particular one of the form $P_{\Phi}\to B \times S^1$.
Consider its
pull-back $P'$ by the projection map $B \times T^2 \to B \times S^1$. This is
still a Hamiltonian bundle. To show that $P$ c-splits, it is sufficient,
by Lemma~\ref{le:funct} $(ii)$, to show that $P'$ c-splits. 
Because $B \times T^2$ may be considered as a smooth compact Hamiltonian
fibration $T^2 \hookrightarrow (B \times T^2) \to B$ with simply connected
base, Lemma~\ref{le:funct2} applies. Thus $P'$ c-splits since any Hamiltonian
bundle over $B$ or over $T^2$ c-splits.  \QED

Finally, we prove  the statements made in~\S\ref{ss:strchar} about the
automorphism groups of Hamiltonian structures.
\MS

\NI
{\bf Proof of Proposition~\ref{prop:hamchar4}}\,\,  

We have  to show that 
the following statements are equivalent for  any $\Phi\in
\Symp_0(P,\pi)$:
\smallskip

$(i)$\,\, $\Phi$ is isotopic to an element of $ \Ham(P, \pi)$;

$(ii)$\,\, 
$\Phi^*(\{a\}) = \{a\}$ for some
Hamiltonian structure $\{a\}$  on $P$;  

$(iii)$
$\Phi^*(\{a\}) = \{a\}$  for all  
Hamiltonian structures $\{a\}$  on $P$. 
\MS

Recall from  Lemma~\ref{le:ham2c} that the
relative homotopy groups $\pi_i(\Symp(M), \Ham(M))$ all vanish for $i >
1$.  Using this together with the fact that $a\in H^2(P)$,  we can reduce to the
case when $B$ is a closed oriented surface.  The statement
 $(i)$ implies $(iii)$  then follows immediately from
Proposition~\ref{prop:automorphisms}. Of course, $(iii)$ implies $(ii)$ so it
remains to show that $(ii)$ implies $(i)$. 

Let us prove this first in the case
where $P\to B$ is trivial, so that $\Phi$ is a map $
B \to  \Symp_0(M, \om)$.  Suppose  that $\Phi^*(a) = a$
for some extension class $a$.   By isotoping $\Phi$ if necessary, we can
suppose that  $\Phi$ takes the base point $b_0$ of $B$ to the identity map. 
Then,  for  each loop $\ga$ in $B$ and any trivialization $T_\ga$,  $$
\begin{array}{cclcl}
0 & =&  f(T_\ga, \Phi^*(a)) - f(T_\ga,a) & = & f(T_\ga\circ\Phi_\ga, a) - f(T_\ga,a)
\\ & = & f(T_\ga, a)\circ \tr_{\Phi_\ga} & = & \om \circ \tr_{\Phi_\ga}\\ &  = &
\Flux(\Phi_\ga) &&
\end{array}$$
where $\Phi_\ga$ is the loop  given by restricting 
$\Phi$ to $\ga$.  Thus the composite
$$
\pi_1(B)\stackrel{\Phi_*}\longrightarrow 
\pi_1(\Symp_0(M))\stackrel{\Flux}\longrightarrow H^1(M,\R) 
$$
must vanish.   Thus the restriction of $\Phi: B\to \Symp_0(M)$ to the
$1$-skeleton of $B$ homotops into $\Ham(M)$. Since the
relative homotopy groups $\pi_i(\Symp(M), \Ham(M))$ all vanish for $i >
1$, this implies that $\phi$ homotops to a map in $\Ham(M)$, as required.

Therefore, it remains to show
that we can  reduce the proof that  $(ii)$ implies $(i)$
to the case when $P\to B$ is trivial.  To this end, 
 isotop $\Phi$ so that it is
the identity map on all fibers $M_b$ over some disc $D\subset B$.  Since
$P\to B$ is trivial over $X = B - D$, we can decompose $P\to B$ into the fiber
connected sum of a trivial bundle $P_B$ over $B$ (where $B$ is thought of as
the space obtained from $X$ by identifying its boundary to a point) and a
nontrivial bundle $Q$ over $S^2 = D/\p D$.   Further, this decomposition is
compatible with $\Phi$, which can be thought of as the fiber sum of some
automorphism $\Phi_B$ of $P_B$ together with the trivial automorphism of
$Q$.  Clearly, this reduces the proof that $(ii)$ implies $(i)$ to the case
$\Phi_B: P_B \to P_B$ on trivial bundles, if we note  that when $\Phi_B$
is the identity over some disc $D \subset B$, the isotopy between $\Phi$ and an element
in $\Ham(P_B)$ can be constructed so that it 
remains equal to the identity over $D$.
 \QED

\section{The cohomology of general symplectic bundles}\label{ss:nonH}

In this section we discuss some consequences
 for general symplectic bundles of
our results on Hamiltonian
bundles.  First, we prove 
Proposition~\ref{prop:p}  that states that  the boundary map $\p$ in the
rational homology Wang sequence of a symplectic bundle over $S^2$ has  
$\p\circ\p =
0$. \MS

\NI
{\bf Proof of Proposition~\ref{prop:p}}\MS

The map $\p\circ\p: H_k(M)\to H_{k+2}(M)$ is given by
$a\mapsto \Psi_*([T^2]\otimes a)$ where
$$
\Psi: T^2\times M\to M: 
(s,t,x)\mapsto \phi_s\phi_t(x).
$$
Let $b(s,t) = \phi_{s+t}^{-1} \phi_s \phi_t.$
The map $(s,t)\mapsto b(s,t)$
factors through $$
f: T^2\; \longrightarrow \;S^2 = T^2/\{s=0\}\cup\{t=0\}.
$$
Let $Z\to M$ represent a k-cycle.
We have a map
$$
T^2\times Z \stackrel{A_1}\to S^1\times S^2\times Z \stackrel{A_2}\to M
$$
given by 
$$
(s,t,z)\mapsto (\phi_{s+t}, f(s,t), z)\mapsto \phi_{s+t} b(s,t)z =
 \phi_s\phi_t(z),
$$
and want to calculate
$$
\int_{T^2\times Z} A_1^* A_2^*(\alpha) = \int_{(A_1)_*[T^2\times Z]}
A_2^*(\alpha)
$$
for some $k+2$-form $\alpha$ on $M$.
But $(A_1)_*[T^2\times Z] \in H_2(S^2)\otimes H_k(Z)$.  (There is no
component in $H_3(S^1\times S^2)\otimes H_{k-1}(Z)$ since $A_1 = id$ on the
$Z$ factor.)  
Now observe that $A_2^*(\alpha)$ vanishes on $H_2(S^2)\otimes H_k(Z)$ by
Theorem~\ref{thm:act}.\QED

The previous lemma is trivially true for any smooth
(not necessarily symplectic)  bundle over
$S^2$ that extends to $\CP^2$.  For the differential $d_2$ in the Leray
cohomology spectral sequence can be written as 
$$
d_2(a) = \p(a)\cup u \in E_2^{2,q-1},
$$
 where $a\in H^q(M) \equiv
E_2^{0,q}$ and $u$ generates $ H^2(\CP^2) \equiv E_2^{2,0}.$ Hence
$$
0 = d_2(d_2(a)) = d_2(\p(a)\cup u) = d_2(\p(a)) \cup u = \p(\p(a)) \otimes
u^2.
$$

\begin{lemma}  If $\pi:P\to B$ is any symplectic bundle over a simply
connected base, then $d_3 \equiv 0$.
\end{lemma}
\proof{}  As in the proof of Lemma~\ref{le:shamchar} we can reduce to the
case when $B$ is a wedge of $S^2$s and $S^3$s.  The differential $d_3$ is
then given by restricting to the bundle over $\vee S^3$.  Since this is
Hamiltonian, $d_3\equiv 0$ by Theorem~\ref{thm:split}. \QED 

\MS

The next lemma describes the Wang differential $\p = \p_\phi$ in the case of
a symplectic loop $\phi$ with nontrivial image in $H_1(M)$.   

\begin{lemma}\label{le:nH} Suppose that $\phi$ is a symplectic loop such
that
  $[\phi_t(x)]\ne 0$ in $H_1(M)$.
Then  $\kerr\p = \im \p$,
where $\p = \p_\phi: H^k(M) \to H^{k-1}(M)$  is the
corresponding Wang differential. 
\end{lemma}
\proof{} Let $\al\in H^1(M)$ be
such that $\al([\phi_t(x)]) = 1$.  So $\p \al= 1$.  Then, for every $\beta\in
\kerr \p$, $\p(\al\cup\beta) = \beta$.  This means that 
$\kerr \p \subset \im \p$ and so $\kerr \p = \im \p$ (using the fact that
$\p\circ\p = 0$.)  Moreover
 the map
$$
\al\cup: H^k(M) \to H^{k+1}(M)
$$
is injective on $\kerr\, \p$ and $H^*(M)$ decomposes as the direct sum
$\kerr\,\p \oplus (\al\cup\kerr \,\p)$.  
\QED

\NI
{\bf Proof of Corollary~\ref{cor:vnH} }

This claims that for a symplectic loop $\phi$,  $\kerr\,\p = \im \,\p$
if and only if $[\phi_t(x)]\ne 0$ in $H_1(M)$.  The above lemma proves the ``if"
statement.  But the ``only if" statement is easy.  Since $1\in H^0( M)$ is in 
$\kerr\,\p$ it must equal $\p(\al)$ for some $\al\in H^1(M)$.  This means that
$\al([\phi_t(x)])\ne 0$ so that $[\phi_t(x)]\ne 0$.\QED 

\begin{rmk}\rm 
The only place that the symplectic condition enters in  the proof of
Lemma~\ref{le:nH} is in the claim that $\p\circ\p = 0$.  Since this is
always true when the loop comes from a circle action,  this lemma
holds for all, not necessarily symplectic, circle actions.    In this case, we can
interpret the result topologically. For the hypothesis  $[\phi_t(x)]\ne 0$ in
$H_1(X)$ implies that the action has no fixed points, so that the quotient
$M/S^1$ is
 an orbifold with cohomology isomorphic to $\kerr\,\p$.  Thus, 
the argument shows that $M$ has the
same cohomology as the product $(M/S^1)\times S^1$. \end{rmk}

\appendix

\section{More on Hamiltonian structures}

Another approach to characterizing a Hamiltonian structure is to define it
in terms of a structure on the fiber that is
preserved  by elements of the Hamiltonian group.  This section developed 
via discussions with Polterovich.

\begin{defn} \rm A {\it marked symplectic manifold} $(M,\om, [L])$
is a pair consisting of a  closed symplectic manifold $(M,\om)$ together with a
marking $[L]$.  Here $L$ is a collection $\{\ell_1,\dots,\ell_k\}$ of
 loops $\ell_i:S^1\to M$ in $M$ that projects to 
 a minimal generating set $\Gg_L = \{ [\ell_1],\dots, [\ell_k]\}$ for
$H_1(M,\Z)/{\rm torsion}$.  A marking $[L]$ is an equivalence class of
generating loops $L$, where  $L\sim L'$ if for each $i$ there is 
an singular integral $2$-chain $c_i$ 
whose boundary modulo torsion is $\ell_i'-\ell_i$  such that
$\int_{c_i}\om = 0$. \end{defn}

The symplectomorphism group acts on the space 
$\Ll$ of markings.  Moreover, it is easy to check that if a symplectomorphism
$\phi$ fixes one marking $[L]$ it fixes them all.
Hence the group
$$
\LHam(M, \om) = \LHam(M, \om, [L]) = \{\phi\in \Symp(M,\om):\phi_*[L] =
[L]\} $$
independent of the choice of $[L]$.  Its identity component is 
$\Ham(M,\om)$.

 There is a forgetful map $[L]\to \Gg_L$ from the space $\Ll$ of markings
to the space of minimal generating sets for the group $H_1(M,\Z)$, and it is not
hard to check that its   fiber is $(\R/\PPp)^k$, where $\PPp$ is the  image of
the period homomorphism  $$ I_{[\om]}: H_2(M, \Z) \to \R.
$$ 
If $\PPp$ is not discrete, there is no nice topology one can
put on $\Ll$.  However,  it has a pseudotopology, i.e. one can specify which maps
of finite polyhedra $X$ into $\Ll$ are continuous, namely:
$f:X\to \Ll$ is continuous if and only if every $x\in X$ has a neighborhood
$U_x$ such that $f:U_x\to \Ll$ lifts to a continuous map into the space of
generating loops $L$.

Here is an alternative formulation of Theorem~\ref{thm:hamchar3} 
in the language of markings.  

\begin{prop}\label{prop:mark}   Fix a marking $[L]$ on $(M, \om)$.
 A Hamiltonian structure on a symplectic bundle $\pi:P\to B$ is 
an isotopy class of symplectomorphisms $(M, \om, [L]) \to (M_{b_0},
\om_{b_0}, [L_0])$
 together with a continuous choice of marking 
$[L_b]$ on each fiber $(M_b,
\om_b)$ 
 that is trivial over the $1$-skeleton $B_1$ of $B$ in the sense that
 there is a symplectic trivialization 
$$
\Phi: \pi^{-1}(B_1) \to B_1\times (M,\om, [L])
$$
that respects the markings on each fiber.
\end{prop}

Here is another way of thinking of a Hamiltonian structure due to 
Polterovich.\footnote{Private communication}  He observed that there is an exact
sequence 
$$
0 \to \R/\PPp\to SH_1(M,\om) \to H_1(M,\Z) \to 0,
$$
where  $SH_1(M,\om)$ is the ``strange homology group" formed by quotienting
the space of  integral $1$-cycles by the image under $d$ of  the space of integral
$2$-chains with zero symplectic area.   The group $\Symp(M,\om)$ acts on 
$SH_1(M,\om)$.  Moreover, if $\phi\in \Symp_0(M)$ and 
$\ta\in SH_1(M)$ projects to $a\in
H_1(M)$, then $\phi_*(\ta) - \ta \in \R/\PPp$ can be thought of as the value of
the class  $\Flux(\phi) \in H^1(M, \R)/\Ga_{\om}$ on $a$.
It is easy to see that  $\LHam(M,\om)$ is  the subgroup
of $\Symp(M,\om)$ that acts trivially on $SH_1(M,\Z)$.  Further, a marking
on  $(M,\om)$ is a pair consisting of a splitting of the above sequence together
with a generating set $\Gg_L$ for $H_1(M,\Z)/{\rm torsion}$.

Given any symplectic bundle $P\to B$ there is  an associated
bundle of abelian groups with fiber $SH_1(M,\om)$. 
 A Hamiltonian structure on $P\to B$ is a
flat connection on this bundle that is trivial over  the $1$-skeleton $B_1$, under
an appropriate equivalence relation.

These ideas can obviously be generalized to bundles that are not trivial over
the $1$-skeleton.  Equivalently, one can consider bundles with
disconnected structural group.  This group could be the whole of
$\LHam(M,\om)$.  One could also restrict to elements
 acting trivially on $H^*(M)$ and/or to those that act trivially on the
groups 
$$
SH_{2k-1}(M,\om) = \frac{\mbox{integral
$(2k-1)$-cycles}}{d(\mbox{$2k$-chains in the kernel of $\om^k$})}.
$$
These generalizations of $SH_1(M,\om)$ are closely connected to
Reznikov's Futaki type characters: see~\cite{Rez}~\S4.  
It is not yet clear what
is the most natural disconnected extension of $\Ham(M,\om)$.

\MS
\NI
{\bf Acknowledgements.}  The authors are deeply grateful to 
Leonid Polterovich
for several discussions concerning many of the questions 
developed in this paper. He pointed out a few
important mistakes at various stages of the work, which was especially helpful
in finding good conditions under which the main results
of this paper hold. We are also grateful to  Steven Boyer, Jarek Kedra and 
Dennis Sullivan for
making  some pertinent remarks, to Igor Belegradek for alerting
us to Meier's work~\cite{Me} and to ETH (Zurich) where part of this work
was undertaken.

\section{Erratum}

\NI
{\bf Note:} This erratum is published in {\it Topology} {\bf 44} (2005), 1301--3. The paper itself is published in
{\it Topology} {\bf 42} (2003), 309--347.
\MS

The statements in this paper that characterize Hamiltonian bundles $(M,\om)\to P\to B$ are not correct when $H_1(B;\Z)$ has torsion.  The affected results are Theorem 1.1, Proposition 1.2 and Lemma 1.4.  The problem is that the proof of Lemma 2.5 works only if $H_1(B;\Z)$ is a free group.  Hence the arguments in this lemma prove the following weaker version of Theorem 1.1.

\begin{prop}  Let $\pi:P\to B$ be a smooth symplectic fiber bundle with fiber $(M,\om)$ that is symplectically trivializable over the $1$-skeleton.  Then the following conditions are equivalent:
\smallskip 

\NI
(i)  there is
a cohomology class $a\in H^2(P,\R)$ that restricts to $[\om]$ on the fiber $M$;
\smallskip 

\NI
(ii) the pullback of $P\to B$ over a suitable finite cover $\Tilde B\to B$ has a Hamiltonian structure.
\end{prop}

The following example (due to Dietmar Salamon) shows that it can be necessary to pass to the finite cover in (ii).  Consider
 the quotient 
$$
P:= \frac{S^2\times\T^2}{\Z_2}
$$
where we think of $S^2\subset\R^3$ as the unit 
sphere and of $\T^2=\R^2/\Z^2$ as the standard torus;
the nontrivial element of $\Z_2$ acts by the involution
$$
(x,y)\mapsto (-x,y+(1/2,0)),\qquad  x\in S^2,\; y\in\T^2.
$$
The closed $2$-form
$
\tau = dy_1\wedge dy_2\in\Om^2(S^2\times\T^2)
$
descends to a closed connection $2$-form on $P$;
its holonomy around each contractible loop in 
$\R P^2$ is the identity and around each noncontractible loop
is the symplectomorphism $(y_1,y_2)\mapsto(y_1+1/2,y_2)$. 
Thus the bundle $\pi: P\to \R P^2$ satisfies the hypothesis
of the above proposition as well as conditions (i) and (ii).  But it does not have a Hamiltonian structure because the classifying map 
$\R P^2\to B\Symp_0(\T^2)$ is not null homotopic, while $\Ham(\T^2)$ 
is contractible.

The mistake in the proof of Lemma 2.5 was the tacit assumption that
the flux class 
$[f(T_\ga,a)] \in H^1(M,\R)/\Ga_{\om}$ vanishes when $\ga\in \pi_1(B)$ has finite order in $H_1(B;\Z)$.
If this condition holds, there is a lifted homomorphism
$\Tilde f_a: \pi_1(B)\to H^1(M,\R)$  such that 
$pr\circ \Tilde f_a([\ga]) = [f(T_\ga,a)]$, where $pr$ denotes the projection, and the proof of Lemma 2.5 goes through.

We claim that the existence of the lift $\Tilde{f}_{a}$ 
does not depend on the choice of extension $a$.
 To see this, choose a symplectic trivialization $T$
 over the $1$-skeleton $B_1$ and a closed connection form $\tau$ in class $a$.  Then $\tau$-parallel translation 
 around a loop $\ga$ in $B$ gives rise to a path $g_t, t\in [0,1],$ 
in $\Symp(M)$ that starts at the identity, and by definition
$$
f(T_\ga,a):= \Flux \bigl(\{g_t\}\bigr) \in H^1(M;\R).
$$
Note that the image $[f(T_\ga,a)]$ of  $f(T_\ga,a)$ in $H^1(F;\R)/\Ga_\om$ is independent of the choice of trivialization. 
The lift $\Tilde f_a$ exists if and only if 
$f(T_\ga,a)$ belongs to  $\Ga_\om$ for all loops $\ga$ with finite order in $H_1(B;\Z)$. 
Given such a loop $\ga$ (which we can assume to be embedded) and a $1$-cycle $\de$ in $M$, denote by
$C(\ga,\de)$ the $2$-cycle in $P$ that equals $\ga\times \de$
under the identification  of $\pi^{-1}(\ga)$ with 
$\ga\times M$ given by $T$.  
Then 
$$
f(T_\ga,a)([\de]) = \int_{C(\ga,\de)} \tau.
$$
But if  $\ga$ has  order $k$ in $H_1(B;\Z)$
the cycle $kC(\ga,\de)$ is homologous to a cycle in the fiber $M$.
Hence
 the integral of $\tau$ over $C(\ga,\de)$ is determined by $[\om]$.
Thus, when $\ga$ is homologically torsion, the class $f(T_\ga,a)$
does not depend on the choice of $a$. Moreover, if $f(T_\ga,a)$ belongs to $\Ga_\om$, one can change the trivialization $T$ over $\ga$ to a trivialization $T'$ such that $f(T_\ga',a)=0$. 
If $f(T_\ga',a)=0$ for all loops $\ga$ that represent a torsion class
 in $H_1(B;\Z)$ we shall say that 
 {\bf the flux of $T'$ vanishes on torsion loops}.

Here are corrected versions of Theorem 1.1 and Proposition 1.2.

\begin{thm}\label{thm:hamchar}  A symplectic bundle
  $\pi: P\to B$ is Hamiltonian
if and only if the
following conditions hold:\begin{itemize}
\item[(i)]  the restriction of $\pi$ to the $1$-skeleton $B_1$ of $B$
has a symplectic trivialization whose flux vanishes on torsion loops, and
\item[(ii)] there is
a cohomology class $a\in H^2(P,\R)$ that restricts to $[\om]$ on the fiber $M$.
\end{itemize}
\end{thm}

\begin{prop}\label{prop:hamchar2}
A  symplectic
bundle  $\pi: P\to B$ is Hamiltonian if and only if the forms $\om_b$ 
on the fibers have  a
closed extension $\tau$ such that the  holonomy of the corresponding connection $\nabla_{\tau}$
around any loop $\ga$ in $B$ 
lies in the identity component $\Symp_0(M)$ of $\Symp(M)$ and moreover lies in $\Ham(M)$ whenever $\ga$ has finite order in $H_1(M;\Z)$. 
\end{prop}

Lemma 1.4 is correct if one understands  it to refer to the corrected version of Theorem 1.1. 
The only other argument that requires comment is the proof of Hamiltonian stability.  
Note first that Corollary 3.2 needs an extra hypothesis
to ensure that the transition functions $\phi_{ij}$ of $P\to B$ preserve the cohomology class $[\om']$ of the perturbed form. 
This hypothesis is satisfied in the Hamiltonian case, since the $\phi_{ij}$ 
may be assumed to be isotopic to the identity.
The next problem is that the proof of Lemma 3.4 uses the incorrect version of Theorem 1.1.  However, 
we can first reduce to the case when $\pi_1(B)=0$ by using Corollary 2.6, and then the two versions of Theorem 1.1 coincide.

\end{document}